\documentclass[showpacs,preprint,superscriptaddress,amsmath,amssymb]{revtex4-1}
\usepackage{amsfonts} 
\usepackage{bbm}
\usepackage{makecell}
\usepackage{dcolumn}% Align table columns on decimal point
\usepackage{bm}% bold math
\usepackage{graphicx}
\usepackage{epstopdf}
\usepackage{booktabs} %Ê¹ÓÃ¼Ó´ÖµÄ±í¸ñÏßÐèÒªÊ¹ÓÃÕâ¸ö°ü
\usepackage{subfig}
\usepackage{makecell}
\usepackage{supertabular}
\usepackage{algorithm}%
\usepackage{algorithmicx}%
\usepackage{amssymb}
\usepackage{amsthm}%Ìá¹©\newtheoremµÈ¶¨Òå¶¨Àí»·¾³µÄÃüÁî
\usepackage{amsmath}
\allowdisplaybreaks[3]
\usepackage{array}
\usepackage{color}
\usepackage{appendix}
\usepackage{bm}
\usepackage{slashed}
\makeatletter
\newcommand{\biggg}{\bBigg@{3}}

\newcommand{\Biggg}{\bBigg@{3.5}}

\newcommand{\bigggg}{\bBigg@{4}}

\newcommand{\Bigggg}{\bBigg@{4.5}}

\makeatother
\usepackage{tabularx}
\usepackage{multirow}

\usepackage{lineno}  %¸øÎÄÕÂ¼ÓÐÐÊý
\usepackage{xcolor}
\usepackage{algorithm}
\usepackage{algpseudocode}
\usepackage{threeparttable}

\begin{document}
	\sloppy
	%\linenumbers       %¸øÎÄÕÂ¼ÓÐÐÊý
	
	\preprint{APS/123-QED}
	
	\title{A Cole-Hopf transformation based fourth-order multiple-relaxation-time lattice Boltzmann model for the coupled Burgers' equations}
	\author{Ying Chen}
	\affiliation{School of Mathematics and Statistics, Huazhong
		University of Science and Technology, Wuhan 430074, China}
	\author{Xi Liu}
	\affiliation{School of Mathematics and Statistics, Huazhong
		University of Science and Technology, Wuhan 430074, China}
	\author{Zhenhua Chai}
	\email[Corresponding author]{(hustczh@hust.edu.cn)}
	\affiliation{School of Mathematics and Statistics, Huazhong
		University of Science and Technology, Wuhan 430074, China}
	\affiliation{Institute of Interdisciplinary Research for Mathematics and Applied Science, Huazhong University of Science and Technology, Wuhan 430074, China}
	\affiliation{Hubei Key Laboratory of Engineering Modeling and
		Scientific Computing, Huazhong University of Science and Technology,
		Wuhan 430074, China}
	
	\author{Baochang Shi}
	\affiliation{School of Mathematics and Statistics, Huazhong
		University of Science and Technology, Wuhan 430074, China}%
	\affiliation{Institute of Interdisciplinary Research for Mathematics and Applied Science, Huazhong University of Science and Technology, Wuhan 430074, China}
	\affiliation{Hubei Key Laboratory of Engineering Modeling and
		Scientific Computing, Huazhong University of Science and Technology,
		Wuhan 430074, China}

	\date{\today}% It is always \today, today,but any date may be explicitly specified
	
	%%%%% Begin Abstract %%%%%%%%%%%
	\begin{abstract}
		In this work, a Cole-Hopf transformation based fourth-order multiple-relaxation-time lattice Boltzmann (MRT-LB) model for   $d$-dimensional coupled Burgers' equations is developed. We first adopt the Cole-Hopf transformation where an intermediate variable $\theta$ is introduced to eliminate 
		the nonlinear convection terms in the Burgers' equations on the velocity $\bm{u}=(u_1,u_2,\ldots,u_d)$. In this case,  a diffusion equation on the variable $\theta$ can be obtained, and particularly, the velocity $\bm{u}$ in the coupled Burgers' equations is determined by the variable $\theta$ and its gradient term $\nabla\theta$. Then we develop the MRT-LB model with the natural moments for the $d$-dimensional transformed diffusion equation and present the corresponding macroscopic finite-difference scheme. At the diffusive scaling, the fourth-order modified equation of the developed MRT-LB model is derived through the Maxwell iteration method.  With the aid of the free parameters in the MRT-LB model, we find that not only the consistent fourth-order modified equation can be obtained, but also the gradient term  $\nabla\theta$ with a fourth-order accuracy can be determined by the non-equilibrium distribution function, this indicates that theoretically,   the   MRT-LB model for $d$-dimensional coupled Burgers' equations can achieve the fourth-order accuracy in space. Finally, some simulations are conducted to test the MRT-LB model, and the numerical results show that the proposed MRT-LB model has a fourth-order convergence rate, which is consistent with our theoretical analysis.
	\end{abstract}
	%%%%% end %%%%%%%%%%%
	\maketitle
	\section{introduction}

	The Burgers' equation, as an important kind of basic partial differential equations (PDEs), can be used to describe some nonlinear physical phenomena, such as turbulence \cite{Konstantin2007}, shock wave \cite{Bogaevsky1990, Reyna1995}, solitons \cite{Buzzicotti2016}, gas dynamics \cite{Lighthill1956}, traffic flow \cite{Nagatani2002}, nonlinear speed of sound \cite{Crighton1979}, to name but a few.   In addition, the coupled Burgers' equations, instead of a single one, are usually adopted to describe some coupled transport phenomena. However, it should be noted that the coupled Burgers' equations are nonlinear, and usually, it is difficult to obtain their analytical solutions. Thus, it is desirable and necessary to develop some numerical methods for the (coupled) Burgers' equation(s). In the past years, two important kinds of numerical approaches have been developed for the (coupled) Burgers' equation(s). The first one is traditional macroscopic numerical methods, such as the finite-difference method \cite{Bahadir1999, Liao2008, Liao2010}, finite-element method \cite{Pandey2009}, and finite-volume method \cite{Gao2016}. The other one is the mesoscopic lattice Boltzmann (LB) method, which is not only a highly efficient second-order kinetic theory-based approach for the fluid flow problems governed by the Navier-Stokes equations \cite{Guo2013, Kruger2017, Wang2019, Succi2008}, but also has been successfully extended to solve the (coupled) Burgers' equation(s) \cite{Boghosian2004, Ma2005, Zene2005, Elton1996, Velivelli2006, Duan2006, Zhang2008, Li2018, Qi2018, Rong2023, Chai2008,Chai2014}. Actually, the existing LB models for the (coupled) Burgers' equation(s) \cite{Boghosian2004, Ma2005, Zene2005, Elton1996, Velivelli2006, Duan2006, Zhang2008, Li2018, Qi2018, Rong2023, Chai2008, Chai2014} can be divided into two main kinds. The first one is  the direct approach where the LB model is developed  to solve the Burgers' equation   \cite{Boghosian2004, Ma2005,	Zene2005, Elton1996, Velivelli2006, Duan2006, Zhang2008, Li2018, Chai2008, Chai2014}, however, this approach may suffer from the numerical instability due to the existence of nonlinear and coupled convection term(s), and is usually limited to the one- or two-dimensional Burgers' equation. The second one  is the indirect approach, which is more universal and stable, and mainly focuses on the $d$-dimensional coupled Burgers' equations \cite{Qi2018, Rong2023}. In the second approach,  the Cole-Hopf transformation \cite{Cole1951}  is used to eliminate the non-linearity and coupling in the $d$-dimensional coupled Burgers' equations, and thus a simple diffusion equation on the variable $\theta$ is obtained. Then one can adopt some numerical methods to solve the transformed diffusion equation  \cite{Yang2019, Mitta2009, Anikonov2015, Zhao2011, Mukundan2019, Zhanlav2016}, and the velocity $\bm{u}$ in the coupled Burgers' equations  is further determined by the variable $\theta$ and   gradient term $\nabla\theta$ which can be calculated conveniently by the first-order moment of  the non-equilibrium distribution function with a second-order accuracy  \cite{Rong2023}.
	%Apart from some traditional macroscopic  numerical methods, such as the multi-step method \cite{Mukundan2019} and the finite-difference method \cite{Zhanlav2016},  Qi et al. \cite{Qi2018} first developed a  Cole-Hopf transformation based  SRT-LB model  for the one-dimensional Burgers' equation. Recently, following the same idea,  Rong et al. \cite{Rong2023} further proposed a   second-order Cole-Hopf transformation based SRT-LB model  for the $d$-dimensional coupled Burgers' equations. %The other one is more stable and focuses on the transformed diffusion equation with variable $\theta$ \cite{Qi2018, Rong2023} based on the Cole-Hopf transformation \cite{Cole1951}, then the variable $\bm{u}$ in the coupled Burgers' equations can be determined from $\theta$ and its gradient term $\nabla\theta$ which can be calculated conveniently with a second-order accuracy by the non-equilibrium distribution function \cite{Rong2023}. %
	
	It should be noted that these works mentioned above all focus on the popular single-relaxation-time LB (SRT-LB) model.  However, the SRT-LB model would be unstable  when the relaxation time is close to 1/2. To overcome this problem,  the advanced multiple-relaxation-time LB (MRT-LB) model  \cite{Ginzburg2002, Chai2020, Chai2023}, which is more general than the SRT-LB model and two-relaxation-time LB (TRT-LB) model \cite{Ginzburg2005}, and can be adopted to achieve better numerical stability and/or accuracy through adjusting some free relaxation parameters \cite{Lallemand2000, Pan2006, Cui2016, Luo2011}. For this reason, Yu et al. \cite{Yu2023} developed a second-order MRT-LB model for the one-dimensional Burgers' equation, and the results are more stable. However, these LB models for the (coupled) Burgers' equation(s)   are only of second-order accuracy, and it is still unclear whether a high-order LB model for the $d$-dimensional ($d\geq 1$) (coupled) Burgers' equation(s) can be obtained. Actually,   we note that some works have been made to develop high-order LB models for the diffusion equations. For example,  Suga \cite{Suga} proposed a fourth-order SRT-LB model with the  D1Q3 lattice structure for the one-dimensional diffusion equation, then Lin et al. \cite{Lin2022} extended this work, and developed a sixth-order MRT-LB model. Recently,  Chen et al. \cite{Chen2023} considered the two-dimensional diffusion equation, and proposed a fourth-order MRT-LB model, where the D2Q5 lattice structure is adopted. However, these high-order LB models are limited to the one- and two-dimensional diffusion equations, and there is no general high-order LB model for the $d$-dimensional diffusion equation. In addition, in the framework of the Cole-Hopf transformation based LB model for the coupled Burgers' equations, we also need to develop a high-order scheme for the  gradient term $\nabla\theta$, besides the variable $\theta$ in the diffusion equation. To this end, in this work, we will propose a general Cole-Hopf transformation based high-order MRT-LB model for  $d$-dimensional ($d\geq 1$) coupled Burgers' equations, then at the diffusive scaling, we  further derive the conditions to ensure that  the MRT-LB model for the $d$-dimensional transformed diffusion equation and the calculation of the velocity $\bm{u}$ in the $d$-dimensional coupled Burgers' equations are of fourth-order accuracy in space.
	% one hand, On the other hand, we will focus on the fourth-order modified equation derived from the macroscopic finite-difference scheme of the MRT-LB model \cite{Be2023, Chen2023-GP} for the transformed diffusion equations, and the fourth-order expansion of the non-equilibrium function \cite{Chai2020}.

	The rest of this paper is organized as follows. In Sec. \ref{Sec-Burger-E}, we first present the $d$-dimensional coupled Burgers' equations and the simple diffusion equation  based on the Cole-Hopf transformation. In Sec. \ref{Sec-LB-model}, we develop a general MRT-LB model for the transformed diffusion equation and present the corresponding macroscopic finite-difference scheme \cite{Be2023}, and derive the fourth-order modified equation and fourth-order expression of the distribution function through the Maxwell iteration method \cite{Chai2020, Chen2023-GP}, %which is used to determine the gradient term $\nabla\theta$
	then the conditions of the fourth-order MRT-LB model for the $d$-dimensional coupled Burgers' equations are given. In Sec. \ref{Sec-Num}, some simulations are carried out to test the accuracy of the developed MRT-LB model, and finally, some conclusions are given in Sec. \ref{Conclusion}.
	\section{The Cole-Hopf transformation for $d$-dimensional coupled Burgers' equations}\label{Sec-Burger-E}
	We now consider the following $d$-dimensional coupled Burgers' equations in the computational domain $\Omega$,
	\begin{align} \label{eq-u}
		\frac{\partial \bm{u}}{\partial t}+\bm{u}\cdot\nabla\bm{u}=\upsilon\nabla^2\bm{u},
	\end{align} 
	which satisfy the following  initial, boundary  and  potential symmetry conditions, 
	\begin{subequations}
		\begin{align}
			&\bm{u}(\bm{x},0) =\bm{\psi}(\bm{x},0), \bm{x}\in\Omega,t>0,\\
			&\bm{u}(\bm{x},t) =\bm{\zeta}(\bm{x},t), \bm{x}\in\partial\Omega,t>0,\\
			&\nabla\bm{u}(\bm{x},t)=(\nabla\bm{u}(\bm{x},t))^T,\bm{x}\in\Omega,t>0, 
		\end{align} 
	\end{subequations}
	where $\partial\Omega$ is the boudary of $\Omega$, $\bm{u}=(u_1,u_2,\ldots,u_d)$ is the velocity to be determined, and is dependent on both space $\bm{x}$[$=(x_1,x_2,\ldots,x_d)\in  \Omega$] and time $t$($>0$). $\bm{\psi}=(\psi_1,\psi_2,\ldots,\psi_d)$ and $\bm{\zeta}=(\zeta_1,\zeta_2,\ldots,\zeta_d)$ are the known $d$-dimensional vector functions, $\upsilon$ is the viscosity coefficient. 
	For the $d$-dimensional coupled Burgers' equations, the Cole-Hopf transformation is given by
	\begin{align} \label{transf}
		\bm{u}(\bm{x},t)=-2\upsilon \frac{1}{\theta(\bm{x},t)}\nabla\theta,\bm{x}\in \Omega,t>0, 
	\end{align}  
	With the help of Eq. (\ref{transf}), the coupled Burgers' equations (\ref{eq-u}) can be reformulated as the following $d$-dimensional diffusion equation with the variable $\theta(\bm{x},t)$,
	\begin{align}\label{eq-theta}
		&\frac{\partial\theta}{\partial t}=\upsilon\nabla^2\theta.
	\end{align}
	
	\section{the fourth-order MRT-LB model for   $d$-dimensional coupled Burgers' equations}\label{Sec-LB-model}		
	In this section, we will first develop a general MRT-LB model for the $d$-dimensional transformed diffusion equation (\ref{eq-theta}). Then based on the previous work \cite{Be2023}, we will perform an analysis on the macroscopic finite-difference scheme of the MRT-LB model, and further derive the consistent fourth-order modified equation through the Maxwell iteration method \cite{Chen2023-GP}. Finally, we will present some details on how to calculate the velocity $\bm{u}$ in the coupled Burgers' equations with a fourth-order accuracy. %Finally, some conditions to ensure that the MRT-LB models for the $d$-dimensional coupled Burgers' equations (\ref{eq-u}) have a fourth-order accuracy will be given.
	\subsection{The MRT-LB model for the $d$-dimensional diffusion equation}
	The evolution of MRT-LB model for the $d$-dimensional  diffusion equation (\ref{eq-theta}) 	can be written as  
	\begin{align}\label{distribution-f}
		f_i(\bm{x}+\bm{c}_i\Delta,t+\Delta t)=f_i(\bm{x},t)-\Big(\bm{M}^{-1}\bm{SM}\Big)_{ik}\big[f_k-f_k^{eq}\big](\bm{x},t),i=1,2,\ldots,q,
	\end{align}
	where $f_i(\bm{x},t)$ and $f_i^{eq}(\bm{x},t)$ are the distribution function and equilibrium distribution at position $\bm{x}$ and time $t$, respectively. $q$ denotes the number of discrete velocities in the D$d$Q$q$ ($q$ discrete velocities in $d$-dimensional space) lattice structure, here we adopt the  D$d$Q$(1+2d^2)$ lattice structure  where the discrete velocity $\bm{c}_i$, the transformation matrix $\bm{M}$ based on the natural moments  and the diagonal relaxation matrix $\bm{S}$ are given by\\
	$d=1$:
	\begin{align}\label{d1q3}
		&	\bm{c}_{x_1}=\big(0,1,-1\big)c, \notag\\ 
		&	\bm{M} =\Big(I_{3},\bm{c}_{x_1}^T,(\bm{c}_{x_1}^{.2})^T\Big)^T,\notag\\ 
		&	\bm{S} ={\rm \makebox{\textbf{diag}}}(s_0,s_1,s_2),\bm{c}=\bm{c}_{x_1},
	\end{align} 
	$d=2$:
	\begin{align}\label{d2q9}
		&	\bm{c}_{x1} =\big(0,1,0,-1,0,1,-1,-1,1\big)c,\notag\\
		&\bm{c}_{x2} =\big(0,0,1,0,-1,1,1,-1,-1\big)c, \notag\\
		&	\bm{M}  =\Big(I_{ 9},\bm{c}_{x1}^T,\bm{c}_{x2}^T,(\bm{c}_{x1}^{.2})^T,(\bm{c}_{x2}^{.2})^T,(\bm{c}_{x1}.\bm{c}_{x2})^T,(\bm{c}_{x1}.\bm{c}_{x2}^{.2})^T,(\bm{c}_{x1}^{.2}\bm{c}_{x2})^T,(\bm{c}_{x1}^{.2}.\bm{c}_{x2}^{.2})^T\Big)^T,\notag\\
		&	\bm{S} ={\rm \makebox{\textbf{diag}}}(s_0,s_1,s_1,s_{21},s_{21},s_{22},1,1,1),\bm{c}=(\bm{c}_{x1}^T,\bm{c}_{x2}^T)^T,
	\end{align} 
	$d=3$:
	\begin{align}\label{d3q19}
		&\bm{c}_{x1} =\big(0,1,0,0,-1,0,0,1,-1,-1,1,1,-1,-1,1,0,0,0,0\big)c,\notag\\
		&\bm{c}_{x2} =\big(0,0,1,0,0,-1,0,1,1,-1,-1,0,0,0,0,1,-1,-1,1\big)c,\notag\\
		&\bm{c}_{x3} =\big(0,0,0,1,0,0,-1,0,0,0,0,1,1,-1,-1,1,1,-1,-1\big)c,\notag\\
		&	\bm{M} =\Big(I_{19},\bm{c}_{x1}^T,\bm{c}_{x2}^T,\bm{c}_{x3}^T,(\bm{c}_{x1}^{.2})^T,(\bm{c}_{x2}^{.2})^T,(\bm{c}_{x3}^{.2})^T,(\bm{c}_{x1}.\bm{c}_{x2})^T,(\bm{c}_{x1}.\bm{c}_{x3})^T,(\bm{c}_{x2}.\bm{c}_{x3})^T,(\bm{c}_{x1}.\bm{c}_{x2}^{.2})^T,\notag\\
		&\qquad(\bm{c}_{x1}^{.2}.\bm{c}_{x2})^T, (\bm{c}_{x1}.\bm{c}_{x3}^{.2})^T,(\bm{c}_{x1}^{.2}.\bm{c}_{x3})^T,(\bm{c}_{x2}.\bm{c}_{x3}^{.2})^T,(\bm{c}_{x2}^{.2}\bm{c}_{x3})^T,(\bm{c}_{x1}^{.2}.\bm{c}_{x2}^{.2})^T,(\bm{c}_{x1}^{.2}.\bm{c}_{x3}^{.2})^T,(\bm{c}_{x2}^{.2}.\bm{c}_{x3}^{.2})^T\Big)^T,\notag\\
		&\bm{S} =\makebox{\textbf{diag}}(s_0,s_1,s_1,s_1,s_{21},s_{21},s_{21},s_{22},s_{22},s_{22},1,1,1,1,1,1,1,1,1),\bm{c}=(\bm{c}_{x1}^T,\bm{c}_{x2}^T,\bm{c}_{x3}^T)^T,
	\end{align}
	$d>3$: 
	\begin{align}\label{ddqq}
		&\bm{c}_{x_1} =\big(Q_1,\underbrace{J_1^{d-1}}_{4(d-1)},O_{\zeta_1}\big)c, \notag\\
		&\bm{c}_{x_2} =\big(Q_2,\underbrace{J_2,O_{4(d-2)}}_{4(d-1)},J_1^{d-2},O_{\zeta_2}\big)c,\notag\\
		&\qquad\qquad\qquad\qquad\qquad\vdots\notag\\
		&\bm{c}_{x_i} =\bigg(Q_i,\underbrace{O_{4(i-1)},J_2,O_{4[d-(i+1)]}}_{4(d-1)},\underbrace{O_{4(i-2)},J_2,O_{4[d-(i+1)]}}_{4(d-2)},\ldots, J_1^{d-i},O_{\eta_i}\bigg)c,\notag\\
		&\qquad\qquad\qquad\qquad\qquad\qquad\qquad\qquad\qquad\vdots\notag\\
		&\bm{c}_{x_k} =\big(Q_k,\underbrace{O_{4(k-1)},J_2,O_{4[d-(k+1)]}}_{4(d-1)},\underbrace{O_{4(k-2)},J_2,O_{4[d-(k+1)]}}_{4(d-2)},\ldots,\underbrace{O_{4(k-i)},J_2,O_{4[d-(k+1)]}}_{4(d-i)},\ldots,J_1^{d-k},O_{\zeta_k}\bigg)c,\notag\\
		&\qquad\qquad\qquad\qquad\qquad\qquad\qquad\qquad\qquad\vdots\notag\\
		&\bm{c}_{x_d} =\bigg(Q_{d},\underbrace{O_{4(d-2)},J_2}_{4(d-1)},\underbrace{O_{4(d-3)},J_2}_{4(d-2)},\ldots,\underbrace{O_{4[d-(i+1)]},J_2}_{4(d-i)},\ldots,\underbrace{O_{4[d-(k+1)]},J_2}_{4(d-k)}, \ldots,O_4,J_2,J_2\bigg)c,\notag\\
		&	\bm{M} =\Big(I_{1+2d^2},\bm{c}_{x_1}^T,\bm{c}_{x_2}^T,\ldots,\bm{c}_{x_d}^T,(\bm{c}_{x_1}^{.2})^T,(\bm{c}_{x_2}^{.2})^T,\ldots,(\bm{c}_{x_d}^{.2})^T,\notag\\
		&\qquad(\bm{c}_{x_1}.\bm{c}_{x_2})^T,(\bm{c}_{x_1}.\bm{c}_{x_3})^T,\ldots,(\bm{c}_{x_1}.\bm{c}_{x_d})^T,\ldots,(\bm{c}_{x_i}.\bm{c}_{x_{i+1}})^T,\ldots,(\bm{c}_{x_i}.\bm{c}_{x_d})^T,\ldots,(\bm{c}_{x_{d-1}}.\bm{c}_{x_d})^T,\notag\\
		&\quad\quad(\bm{c}_{x_1}.\bm{c}_{x_2}^{.2})^T,(\bm{c}_{x_1}.\bm{c}_{x_3}^{.2})^T,\ldots,(\bm{c}_{x_1}.\bm{c}_{x_d}^{.2})^T,\ldots,(\bm{c}_{x_i}.\bm{c}_{x_{i+1}}^{.2})^T,\ldots,(\bm{c}_{x_i}.\bm{c}_{x_d}^{.2})^T,\ldots,(\bm{c}_{x_{d-1}}.\bm{c}_{x_d}^{.2})^T,\notag\\
		&\quad\quad(\bm{c}_{x_1}^{.2}.\bm{c}_{x_2})^T,(\bm{c}_{x_1}^{.2}.\bm{c}_{x_3})^T,\ldots,(\bm{c}_{x_1}^{.2}.\bm{c}_{x_d})^T,\ldots,(\bm{c}_{x_i}^{.2}.\bm{c}_{x_{i+1}})^T,\ldots,(\bm{c}_{x_i}^{.2}.\bm{c}_{x_d})^T,\ldots,(\bm{c}_{x_{d-1}}^{.2}.\bm{c}_{x_d})^T,\notag\\
		&\quad\quad(\bm{c}_{x_1}^{.2}.\bm{c}_{x_2}^{.2})^T,(\bm{c}_{x_1}^{.2}.\bm{c}_{x_3}^{.2})^T,\ldots,(\bm{c}_{x_1}^{.2}.\bm{c}_{x_d}^{.2})^T,\ldots,(\bm{c}_{x_i}^{.2}.\bm{c}_{x_{i+1}}^{.2})^T,\ldots,(\bm{c}_{x_i}^{.2}.\bm{c}_{x_d}^{.2})^T,\ldots,(\bm{c}_{x_{d-1}}^{.2}.\bm{c}_{x_d}^{.2})^T\Big),\notag\\
		&\bm{S} =\makebox{\textbf{diag}} (s_0,s_1I_d,s_{21}I_d,s_{22}I_{d(d-1)/2},I_{3d(d-1)/2}),\bm{c}=(\bm{c}_{x_1}^T,\bm{c}_{x_2}^T,\ldots,\bm{c}_{x_d}^T)^T,
	\end{align}  
	with
	\begin{subequations}
		\begin{align}
			&Q_p=(0,0,\ldots,\underbrace{1}_{{p+1}_{th}},0,\ldots,0,\ldots,\underbrace{-1}_{(d+p+1)_{th}},0,\ldots,0)\in R^{(2d+1)\times 1},\\
			&J_i^p=J_i\otimes I_p=(\underbrace{J_i,J_i,\ldots,J_i}_{p})\in R^{4p\times1}, J_1=(1,-1,-1,1),J_2=(1,1,-1,-1),p\in\mathcal{N},\\
			&O_p=\big(0,0,\ldots,0)^T\in R^{p\times 1},I_{p}=\big(1,1,\ldots,1)^T\in R^{p\times 1}, \\
			&\eta_i=2d(d-1)-2i(2d-1-i),\\ 
			&\bm{c}_{x_1}^{.p_1}.\bm{c}_{x_2}^{.p_1}.\ldots.\bm{c}_{x_d}^{.p_d}=\Big(\prod_{i=1}^d \bm{c}_{xi}(1)^{p_i},\prod_{i=1}^d\bm{c}_{xi}(2)^{p_i},\ldots,\prod_{i=1}^d \bm{c}_{xi}(q)^{p_i}\Big),p_i\in\mathcal{N},i=1,2,\ldots,d,
		\end{align}
	\end{subequations}
	where $c=\Delta x/\Delta t$ is the lattice speed with the lattice spacing $\Delta x$ and time space $\Delta t$,  $\bm{c}_{xi}(i)$ represents the $i_{th}$ element of vector $\bm{c}_{xi}$,  and the relaxation parameters $s_0$, $s_1$, $s_{21}$  and $s_{22}$ are located in the range of  $(0,2)$. In order to recover the diffusion equation (\ref{eq-theta}) correctly from the MRT-LB model (\ref{distribution-f}), the equilibrium distribution function $f_i^{eq}$ should be designed as
	\begin{align}\label{feq}
		f_i^{eq}(\bm{x},t)=w_i\theta(\bm{x},t),
	\end{align}
	which satisfies the following moment conditions \cite{Chai2020},
	\begin{subequations}\label{conseved-moment}
		\begin{align}
			&\sum _i f_i^{eq}(\bm{x},t)=\sum _i f_i(\bm{x},t)=\theta(\bm{x},t),\\
			&\sum _i\bm{c}_i f_i^{eq}(\bm{x},t)=\bm{0},\\
			&\sum _i \bm{c}_i\bm{c}_if_i^{eq}(\bm{x},t)=c_s^2\theta(\bm{x},t)\bm{I}, 
		\end{align}
	\end{subequations}
	where $w_i$ is the weight coefficient, $c_s^2=[2w_1+4(d-1)w_{1+2d}]c^2$ in the D$d$Q$(1+2d^2)$ lattice structure. For simplicity, we only consider the following commonly used conditions that the weight coefficients should satisfy 
	\begin{align}\label{d1-3-weight}
		&w_1 =w_2=\ldots=w_{2d},w_{2d+1}=w_{2d+2}=\ldots=w_{2d^2},w_0=1-2dw_1-2d(d-1)w_{2d+1}.
	\end{align}   
	\subsection{The macroscopic finite-difference scheme and modified equation of the MRT-LB model }\label{PartB}
	In this part, we will first present the macroscopic finite-difference scheme  of the  developed MRT-LB model (\ref{distribution-f}) with the D$d$Q$(1+2d^2)$ lattice structure \cite{Be2023}, and then derive the consistent fourth-order modified equation  through the Maxwell iteration method \cite{Chen2023-GP}. 
	
	Following the previous work \cite{Be2023}, one can obtain the following macroscopic finite-difference scheme of the MRT-LB model (\ref{distribution-f}) for the $d$-dimensional diffusion equation (\ref{eq-theta}),
	\begin{equation}\label{det-adj}
		\det{\Big(T^1_{\Delta t}\bm{I}-\bm{A}\Big)}\bm{m}={\rm adj}\Big(T^1_{\Delta t}\bm{I}-\bm{A}\Big)\bm{B}\bm{m}^{eq}.
	\end{equation} 
	where $\det(\cdot)$ and adj($\cdot$) denote the determinant and adjufate matrix, respectively. $T_{\Delta t}^1[h(\bm{x},t)]:=h(\bm{x},t+\Delta t)$  is the time shift operator  with $h(\bm{x},t)$ representing an arbitrary function dependent on the space $\bm{x}$ and time $t$ \cite{Be2023}, the matrices $\bm{A}$ and $\bm{B}$, the moments  $\bm{m}$ and $\bm{m}^{eq}$ are given by
	\begin{align} 
		&\bm{A}=\bm{W}\Big(\bm{I}-\bm{S}\Big),\bm{B}=\bm{WS}, \notag\\
		&\bm{m}=\bm{Mf},\bm{m}^{eq}=\bm{Mf}^{eq}, 
	\end{align}
	with 
	\begin{align}   \label{FD-scheme-quan}
		&\bm{W}=\bm{M}\bm{T}\bm{M}^{-1},\bm{T}=\makebox{\textbf{diag}} \Big(T^{-\bm{c}_{1}}_{\Delta t},T^{-\bm{c}_{2}}_{\Delta t},\ldots,T^{-\bm{c}_{1+2d^2}}_{\Delta t}\Big), \notag\\
		&\bm{f}=\big(f_1,f_2,\ldots,f_{1+2d^2}\big),\bm{f}^{eq}=\big(f_1^{eq},f_2^{eq},\ldots,f_q^{e{1+2d^2}}\big),
	\end{align}
	
	Then due to the equivalence between the MRT-LB model (\ref{distribution-f}) and  macroscopic finite-difference scheme (\ref{det-adj}), we adopt the Maxwell iteration method to derive the  fourth-order modified equation of the finite-difference scheme (\ref{det-adj})   rather than the  second-order modified equation given by some commonly used asymptotic analysis methods  \cite{Chai2020}. In particular, the  fourth-order modified equation of the finite-difference scheme (\ref{det-adj}) is also equivalent to the fourth-order modified equation obtained from the following equation \cite{Chen2023-GP},
	\begin{equation}\label{Xi}
		\begin{aligned}
			\bm{\Xi}:=\bm{m}-\bigg(\sum_{k=0}^{+\infty}\bm{\Gamma}^k\bigg)\bm{m}^{eq}=\bm{0}, 
		\end{aligned}
	\end{equation}
	with 
	\begin{align}\label{Gamma}
		\bm{\Gamma}=-\bm{S}^{-1}\Big(T_{\Delta t}^1\bm{W}^{-1}-\bm{I}\Big).
	\end{align}
	For the time shift operator $T^1_{\Delta t}$ and matrix $\bm{W}^{-1}$ in Eq. (\ref{Gamma}), we can rewrite them as the   two series expansion expressions,
	\begin{subequations}\label{Time-op-W-1}
		\begin{align}
			T^1_{\Delta t}=&\sum_{k=0}^{+\infty}\frac{\eta^k\Delta x^{2k}\partial_t^k}{k!}=1+\eta\Delta x^2\partial_t+\frac{\eta^2\Delta x^4}{2}\partial_{t}^2+O(\Delta x^6),\\
			\bm{W}^{-1}&=\exp(\Delta x\bm{\mathcal{W}})= \sum_{k=0}^{+\infty}\bigg[\sum_{k=1}^{+\infty}\frac{\Delta x^k\bm{\mathcal{W}}^k}{k!}\bigg]\notag\\
			&=\bm{I}+\Delta x\bm{\mathcal{W}}+\frac{\Delta x^2}{2}\bm{\mathcal{W}}^2+\frac{\Delta x^3}{6}\bm{\mathcal{W}}^3+\frac{\Delta x^4}{24}\bm{\mathcal{W}}^4+O(\Delta x^5),
		\end{align}
	\end{subequations}
	with 
	\begin{align}\label{W}
		&\bm{\mathcal{W}}=\bm{M}{\rm \makebox{\textbf{diag}}}\Big(\bm{e}_1\cdot\nabla,\bm{e}_2\cdot\nabla,\ldots,\bm{e}_{1+2d^2}\cdot\nabla\Big)\bm{M}^{-1},\bm{e}_i=\frac{\bm{c}_i}{c},i=1,2,\ldots,1+2d^2,
	\end{align}
	where the diffusive scaling, i.e., 
	$\Delta t=\eta \Delta x^2$ ($\eta\in \mathcal{R}$), has been used. In addition, with the aid of Eq. (\ref{Time-op-W-1}), we also have
	\begin{subequations}\label{gamma-1-4}
		\begin{align}
			\bm{\Gamma}^1=&-\bm{S}^{-1}\bigg[\Delta x\bm{\mathcal{W}}+\frac{\Delta x^2}{2}\bm{\mathcal{W}}^2+\frac{\Delta x^3}{6}\bm{\mathcal{W}}^3+\frac{\Delta x^4}{24}\bm{\mathcal{W}}^4\notag\\
			&\quad+\eta\Delta x^2\partial_t\Big(\Delta x\bm{\mathcal{W}}+\frac{\Delta x^2}{2}\bm{\mathcal{W}}^2\Big)+\frac{\eta\Delta x^4\partial_{t}^2}{2}\bm{I}\bigg]+O(\Delta x^5),\\
			\bm{\Gamma}^2=&\Bigg(-\bm{S}^{-1}\bigg[\Delta x\bm{\mathcal{W}}+\frac{\Delta x^2}{2}\bm{\mathcal{W}}^2+\frac{\Delta x^3}{6}\bm{\mathcal{W}}^3 +\eta\Delta x^2\partial_{t}\Big( \bm{I}+\Delta x\bm{\mathcal{W}}\Big) \bigg]\Bigg)^2+O(\Delta x^5)\notag\\
			=&\Delta x^2\bm{S}^{-1}\bm{\mathcal{W}}\bm{S}^{-1}\bm{\mathcal{W}}+\Delta x^3\bm{S}^{-1}\bm{\mathcal{W}}\bm{S}^{-1}\Big(\frac{\bm{\mathcal{W}}^2}{2}+\eta\partial_t\bm{I}\Big)+\Delta x^3\bm{S}^{-1}\Big(\frac{\bm{\mathcal{W}}^2}{2}+\eta\partial_t\bm{I}\Big)\bm{S}^{-1}\bm{\mathcal{W}}\notag\\
			&\qquad+\Delta x^4\bm{S}^{-1}\bm{\mathcal{W}}\bm{S}^{-1}\Big(\frac{\bm{\mathcal{W}}^3}{6}+\eta\bm{\mathcal{W}}\partial_t\Big)\notag\\
			&\qquad+\Delta x^4\bm{S}^{-1}\Big(\frac{\bm{\mathcal{W}}^3}{6}+\eta\bm{\mathcal{W}}\partial_t\Big)\bm{S}^{-1}\Big(\frac{\bm{\mathcal{W}}^3}{6}+\eta\bm{\mathcal{W}}\partial_t\Big)\notag\\
			&\qquad+\Delta x^4\bm{S}^{-1}\Big(\frac{\bm{\mathcal{W}}^3}{6}+\eta\bm{\mathcal{W}}\partial_t\Big)\bm{S}^{-1}\bm{\mathcal{W}}+O(\Delta x^5),\\
			\bm{\Gamma}^3=&\Bigg(-\bm{S}^{-1}\bigg[\Delta x\bm{\mathcal{W}}+\frac{\Delta x^2}{2}\bm{\mathcal{W}}^2  +\eta\Delta x^2\partial_t\bm{I}  \bigg]\Bigg)^3+O(\Delta x^5)\notag\\
			=&-\Delta x^3\bm{S}^{-1}\bm{\mathcal{W}}\bm{S}^{-1}\bm{\mathcal{W}}\bm{S}^{-1}\bm{\mathcal{W}}-\Delta x^4\bm{S}^{-1}\bm{\mathcal{W}}\bigg[\Big(\bm{S}^{-1}\bm{\mathcal{W}}\bm{S}^{-1}\Big(\frac{\bm{\mathcal{W}}^2}{2}+\eta\partial_t\bm{I}\Big)\\
			&\qquad+\Big(\frac{\bm{\mathcal{W}}^2}{2}+\eta\partial_t\bm{I}\Big)\bm{S}^{-1}\bm{\mathcal{W}}\bigg]+O(\Delta x^5),\notag\\
			\bm{\Gamma}^4=&\Delta x^4\bm{S}^{-1}\bm{\mathcal{W}}\bm{S}^{-1}\bm{\mathcal{W}}\bm{S}^{-1}\bm{\mathcal{W}}\bm{S}^{-1}\bm{\mathcal{W}}+O(\Delta x^5).
		\end{align}
	\end{subequations}
	Here it should be noted that we only expand the time shift operator $T^1_{\Delta t}$,  $\bm{W}^{-1}$  and $\bm{\Gamma}^k$ $(k=1,2,3,4)$ up to $O(\Delta x^5)$, which is sufficient to derive the fourth-order modified equation of the finite-difference scheme (\ref{det-adj}), i.e., the following truncation equations up to $O(\Delta x^5)$ of Eq. (\ref{Xi}),
	\begin{equation}\label{fourth-Xi}
		\begin{aligned}
			\bm{m}-\bm{m}^{eq}-\bigg(\sum_{k=1}^{4}\bm{\Gamma}^k\Bigg)\bm{m}^{eq}=O(\Delta x^5) .
		\end{aligned}
	\end{equation} 
	
	Based on the results presented above, we can give the detailed fourth-order modified equations of the MRT-LB models (\ref{distribution-f}) with the D$d$Q$(1+2d^2)$ lattice structure for the $d$-dimensional diffusion equation. 
	
	With some algebraic and symbolic manipulations,  it is easy to calculate the moments $\bm{m}^{eq}$ and $\bm{m}$, $\bm{\mathcal{W}}$  and $\bm{\Gamma}^k$ $(k=1,2,3,4)$ from the equilibrium distribution function $\bm{f}^{eq}$, moment conditions, lattice structure, transform matrix $\bm{M}$ and  relaxation matrix $\bm{S}$ of the MRT-LB model (\ref{distribution-f}). Then selecting the first row of Eq. (\ref{fourth-Xi}) corresponding to the conservative moment $\theta$ yields the following fourth-order modified equations,\\
	$d=1$:\begin{align}\label{Me-1d}
		\frac{\partial \theta}{\partial t}=2w_1\Big(\frac{1}{s_1}-\frac{1}{2}\Big)\frac{\Delta x^2}{\Delta t}\frac{\partial^2\theta}{\partial x^2}+\Delta x^2\frac{w_1(s_1-2)R_1}{12s_1^3s_2}\frac{\partial^4\theta}{\partial x^4}+O(\Delta x^4),
	\end{align}
	with 
	\begin{align}\label{Con-1d}
		R_1&=6\big[s_{2}(s_1-1)+ s_1^2+s_{2}-2s_1 \big] - s_{2}s_1^2  +6w_1\big[2s_1(2-s_1) + s_{2}(4+s_1^2-6s_1) \big] , 
	\end{align} 
	where
	\begin{subequations}
		\begin{align} 
			&\frac{\partial ^2\theta}{\partial t^2}=\bigg[2w_1\Big(\frac{1}{s_1}-\frac{1}{2}\Big)\frac{\Delta x^2}{\Delta t}\bigg]^2\frac{\partial ^4\theta}{\partial x^4}+O(\Delta x^2),\\
			&\frac{\partial ^3\theta}{\partial x^2\partial t}=\bigg[2w_1\Big(\frac{1}{s_1}-\frac{1}{2}\Big)\frac{\Delta x^2}{\Delta t}\bigg]\frac{\partial ^4\theta}{\partial x^4}+O(\Delta x^2),
		\end{align}
	\end{subequations}
	$d>1$:\begin{align}\label{Me-d}
		\frac{\partial \theta}{\partial t}=& \Big(\frac{1}{s_1}-\frac{1}{2}\Big)c_s^2\Delta t\Big(\frac{\partial^2\theta}{\partial x_1^2}+\frac{\partial^2\theta}{\partial x_2^2}+\ldots+\frac{\partial^2\theta}{\partial x_d^2}\Big)\notag\\
		& +\Delta x^2\frac{\big[w_1+2(d-1)w_{1+2d}\big](s_1-2)R_{d1}}{12s_1^3s_{21}}\Big(\frac{\partial^4\theta}{\partial x_1^4}+\frac{\partial^4\theta}{\partial x_2^4}+\ldots+\frac{\partial^4\theta}{\partial x_d^4}\Big)\notag\\
		&-\Delta x^2\frac{R_{d2}}{s_1^3s_{21}s_{22}}\Big(\frac{\partial^4\theta}{\partial x_1^2\partial x_2^2}+\frac{\partial^4\theta}{\partial x_1^2\partial x_3^2}+\ldots\frac{\partial^4\theta}{\partial x_1^2\partial x_d^2}+\ldots\frac{\partial^4\theta}{\partial x_{d-1}^2\partial x_{d}^2}\Big)+O(\Delta x^4),
	\end{align} 
	with   
	\begin{subequations}\label{Con-d}
		\begin{align}
			R_{d1}=&6\big[s_{21}(s_1-1)+ s_1^2+s_{21}-2s_1 \big] - s_{21}s_1^2 \notag\\
			&  +6\big[w_1+2(d-1)w_{1+2d}\big]\big[2s_1(2-s_1) + s_{21}(4+s_1^2-6s_1) \big]  ,	\\
			R_{d2} =&s_1^2w_{1+2d}\big[2(s_1-2)(2s_{21} + s_{22}) + s_{21}s_{22}(8-3s_1)  \big]-8s_1s_{21}s_{22}(w_1 + 2(d-1)w_{1+2d})^2\notag\\
			& +s_{22}\Big(w_1^2 + 4(d-1)w_{1+2d}\big[w_1+(d-1)w_{1+2d}\big]\Big)\times\notag\\
			& \Big(8\big[s_1  + s_{21}- s_1^2 + (s_1-1)s_1s_{21}\big]+(2-s_{21})s_1^3\Big),
		\end{align}
	\end{subequations}
	where
	\begin{subequations}
		\begin{align}
			&\frac{\partial ^2\theta}{\partial t^2}=\bigg[\Big(\frac{1}{s_1}-\frac{1}{2}\Big)c_s^2\Delta t\bigg]^2\Big(\frac{\partial ^4\theta}{\partial x_1^4}+\frac{\partial ^4\theta}{\partial x_2^4}+\ldots+\frac{\partial ^4\theta}{\partial x_d^4} \Big)+O(\Delta x^2),\\
			&\frac{\partial ^3\theta}{\partial x_i^2\partial t}=\Big(\frac{1}{s_1}-\frac{1}{2}\Big)c_s^2\Delta t \Big(\frac{\partial ^4\theta}{\partial x_i^2\partial x_{1}^2}+\frac{\partial ^4\theta}{\partial x_i^2\partial x_{2}^2}+\ldots+\frac{\partial ^4\theta}{\partial x_i^2\partial x_{i-1}^2}\notag\\
			&\qquad\qquad\qquad\qquad+\frac{\partial ^4\theta}{\partial x_{i}^4}+\frac{\partial ^4\theta}{\partial x_i^2\partial x_{i+1}^2}+\ldots+\frac{\partial ^4\theta}{\partial x_i^2\partial x_d^2}\Big)+O(\Delta x^2), i=1,2,\ldots,d,
		\end{align}
	\end{subequations}
	have been used.
	
	We now give a remark on Eqs. (\ref{Me-1d}), (\ref{Con-1d}), (\ref{Me-d}) and (\ref{Con-d}).\\
	\textbf{Remark 1.} We would like to point out that the number of the free parameters in the MRT-LB model (\ref{distribution-f}) is more than that of the conditions to ensure that the MRT-LB model (\ref{distribution-f}) for the diffusion equation (\ref{eq-theta}) is fourth-order accurate. For instance, for the one-dimensional diffusion equation, the number of the free parameters in the MRT-LB model   is three (the weight coefficient $w_1$, the relaxation parameters $s_1$ and $s_2$), while the number of the fourth-order conditions is two ($\upsilon=2w_1(1/s_1-1/2)\Delta x^2/\Delta t$ and $R_1=0$). In addition,  the consistent fourth-order modified equation of the MRT-LB model  (or the fourth-order MRT-LB model) for the $d$-dimensional diffusion equation can be obtained once the free parameter in the MRT-LB model (\ref{distribution-f}) satisfy the fourth-order conditions, i.e., $d=1$: $R_1=0$; $d>1$: $R_{d1}=0$ and  $R_{d2}=0$. 
	\subsection{The calculation of the velocity $\bm{u}$ in the coupled Burgers' equation}
	The  aim of this work is to obtain a fourth-order MRT-LB model for the $d$-dimensional coupled Burgers' equations (\ref{eq-u}), this means that it is necessary to determine  the velocity $\bm{u}$ with a fourth-order accuracy in space. Actually, according to the Maxwell iteration method \cite{Chai2020}, the distribution function $\bm{f}$ can be expressed as  
	\begin{align}\label{f-expression}
		\bm{f}=&\bm{f}^{eq}-\eta\Delta x^2\bm{\Lambda}^{-1}\bm{D}\bm{f}^{eq}+\eta^2\Delta x^4\bm{\Lambda}^{-1}\bm{D}\Bigg[\bigg(\bm{I}-\frac{\bm{\Lambda}}{2}\bigg)\bm{\Lambda}^{-1}\bm{Df}^{eq}\Bigg]\notag\\
		&-\eta^3\Delta x^6\bm{\Lambda}^{-1}\bm{D}\Bigg[\frac{\bm{\Lambda}^{-1}\bm{D}^2}{2}+\frac{\bm{D}\bm{\Lambda}^{-1}\bm{D}}{2}-\frac{\bm{D}^2}{6}-\bm{\Lambda}^{-1}\bm{D}\bm{\Lambda}^{-1}\bm{D}\Bigg]\bm{f}^{eq}+...,
	\end{align} 
	with  
	\begin{subequations}
		\begin{align}
			&\bm{\Lambda}=\bm{M}\bm{S}^{-1}\bm{M}^{-1},\\
			&\bm{D}={\rm \makebox{\textbf{diag}}}\big(\partial_t+\bm{c}_1\cdot\nabla,\partial_t+\bm{c}_2\cdot\nabla,\ldots,\partial_t+\bm{c}_{1+2d^2}\cdot\nabla\big). 
		\end{align}
	\end{subequations}
	Here we would like to point out that there will still be some free parameters in the MRT-LB model after obtaining the consistent fourth-order modified equation of the transformed diffusion equation (see Remark 1 for details). Actually, from Eq. (\ref{f-expression}) one can see that although there are mixed partial terms above the second-order,  they can be eliminated with the aid of the free parameters.  For clarity, we give a theorem below.\\
	\noindent\textbf{Theorem 1.} For the $d$-dimensional coupled Burgers' equations (\ref{eq-u}),  once the following conditions are satisfied,\\
	$d=1$:\begin{align}\label{nabla-theta-con-d1}
		{\rm Con}_1= 3s_1\big[s_{2} + s_1-2\big]  - s_1^2 s_{2} +6w_1 (2-s_1)\big[s_1+s_{2}(1-s_1) \big]  =0,	
	\end{align}
	$d>1$: 
	\begin{align} \label{nabla-theta-con-d2}
		\left\{\begin{aligned}
			{\rm Con}_{d1}=&  3s_1\big[s_{21} + s_1-2\big]  - s_1^2 s_{21} +6\big[w_1+2(d-1)w_{1+2d}\big] (2-s_1)\big[s_1+s_{21}(1-s_1) \big]  =0,\\
			{\rm Con}_{d2}=&-\Big(s_{22}w_1^2+4s_{22}\big[w_1+(d-1)w_{1+2d}\big]w_{1+2d}(d - 1)\Big)(s_1 - 2)(s_1 + s_{21} -s_1s_{21})\\
			&-s_1^2w_{1+2d}(2s_{21} + s_{22}- 2s_{21}s_{22})=0,
		\end{aligned}\right.
	\end{align} 
	the distribution function $\bm{f}$ in the MRT-LB model (\ref{distribution-f}) at the diffusive scaling can be simplified as  
	\begin{equation}\label{f-feq-o5}
		\bm{f}=\bm{f}^{eq}-\eta\Delta x^2\bm{\Lambda}^{-1}\bm{\overline{D}}\bm{f}^{eq} +O(\Delta x^5),
	\end{equation}
	where $\bm{\overline{D}}=$\textbf{diag}$\big(\bm{c}_1\cdot\nabla,\bm{c}_2\cdot\nabla,\ldots,\bm{c}_{1+2d^2}\cdot\nabla\big)$.
	In addition  from above Eq. (\ref{f-feq-o5}) and   the moment condition (\ref{conseved-moment}), we can obtain
	\begin{align}
		\sum_{k=1}^{1+2d^2}\bm{c}_{k}(f_k-f_k^{eq}):=\sum_{k=1}^{1+2d^2}\bm{c}_{k}f^{ne}_k=-\frac{1}{s_1}\Delta t\Big( \nabla\cdot\sum_{k=1}^{1+2d^2}\bm{c}_k\bm{c}_kf_k^{eq}\Big)=-\frac{1}{s_1}\Delta tc_s^2\nabla\theta,
	\end{align}  which means that the gradient term $\nabla\theta$ can be determined by
	\begin{align}\label{nabla-theta}
		\nabla\theta=-\sum_{k=1}^{1+2d^2}\bm{c}_{k}f^{ne}_ks_1\Delta t c_s^2,
	\end{align}
	which is fourth-order accurate in space. Then  according to Eqs. (\ref{transf}) and (\ref{conseved-moment}), the velocity $\bm{u}$ in the $d$-dimensional coupled Burgers' equations (\ref{eq-u}) can be calculated with the fourth-order accuracy  by the following formula,
	\begin{align}\label{ui-compute}
		u_i=-2\upsilon\frac{\partial_{x_i}\theta}{\theta(\bm{x},t)}=\big(2-s_1\big)\frac{\sum_{k=1}^{1+2d^2} c_{ki}f^{ne}_k(\bm{x},t)}{\sum_{k=1}^{1+2d^2}f_k(\bm{x},t)},i=1,2,\ldots,d. 
	\end{align}
	\subsection{The conditions of the fourth-order MRT-LB model}
	Now we give a theorem to present the conditions that ensure the MRT-LB model for the $d$-dimensional coupled Burgers' equations, i.e., Eqs. (\ref{distribution-f}) and (\ref{ui-compute}), to have a fourth-order accuracy in space. \\
	\noindent\textbf{Theorem 2.} The Cole-Hopf transformation based MRT-LB model  [Eqs. (\ref{distribution-f}) and (\ref{ui-compute})] for the $d$-dimensional coupled Burgers' equations (\ref{eq-u}) can be fourth-order accurate once the weight coefficients and the relaxation parameters are given by
	\begin{align}\label{Bur-Con}
		\begin{matrix}   
			d=1:\:\left\{\begin{aligned}
				&\varepsilon=2w_1\xi,\\
				&R_1=0,\\
				&{\rm Con}_1=0,
			\end{aligned}\right.  &	
			d>1:\:\left\{\begin{aligned}
				&\varepsilon=\big[2w_1+4(d-1)w_{1+2d}\big]\xi, \\
				&R_{d1}=0,\\
				&R_{d2}=0,\\
				&{\rm Con}_{d1}=0,\\
				&{\rm Con}_{d2}=0,
			\end{aligned}\right.
		\end{matrix}
	\end{align}
	where $\xi=1/s_1-1/2$ and $\varepsilon:=\upsilon \Delta t/\Delta x^2$. From Eq. (\ref{Bur-Con}), one can determine  the weight coefficients and relaxation parameters as\\
	\begin{align}\label{Bur-Sol}
		\begin{matrix}   
			d=1:\:\left\{\begin{aligned}
				&s_1=2/(6\varepsilon+1),\\
				&s_2=24\varepsilon/(6\varepsilon+1)^2,\\
				&w_1=1/6,
			\end{aligned}\right.  &	
			d>1:\:\left\{\begin{aligned}
				&s_1=2/(6\varepsilon+1),\\
				&s_{21}=24\varepsilon/(6\varepsilon+1)^2,\\
				&s_{22}=4/(6\varepsilon+3),\\
				&w_1=1/6-(d-1)\varepsilon/3,\\
				&w_{1+2d}=\varepsilon/6.
			\end{aligned}\right. 
		\end{matrix}
	\end{align}
	Additionally,  the gradient term $\nabla\theta$  in the diffusion equation (\ref{eq-theta}) can also be calculated locally by the non-equilibrium distribution function  [see  Eq. (\ref{nabla-theta})] with a fourth-order accuracy in space. %Then the variable $\bm{u}$ can be solved very conveniently by $\theta$ and $\nabla\theta$ as shown in Eq. (\ref{ui-compute}).

	\section{NUMERICAL RESULTS AND DISCUSSION}\label{Sec-Num}
	In this section, we will conduct some numerical experiments with the fourth-order MRT-LB model  where the weight coefficient and the relaxation parameters are determined from Eqs. (\ref{d1-3-weight}) and (\ref{Bur-Sol}).  To evaluate the difference between the analytical  and numerical results, the   root-mean-square error (RMSE) is adopted,
	\begin{align}
		{\rm RMSE}=\sqrt{\frac{\sum_{i=1}^d\sum_{j_i=1}^{N_{x_i}}\big[\psi(j_1\Delta x,\ldots,j_{d}\Delta x ,n\Delta t)-\psi^{\star}(j_1\Delta x,\ldots,j_{d}\Delta x ,n\Delta t)\big]^2}{\prod_{i=1}^dN_{x_i}}}
	\end{align}
	where $N_{x_i}$ ($i=1,2,\ldots,d$) is the number of grid points, $\psi$ and $\psi^{\star}$ are the numerical and analytical solutions. In order to estimate the convergence rate (CR) of the MRT-LB model, the following formula is used \cite{Kruger2017},
	\begin{align}
		{\rm CR}=\frac{\log\big({\rm RMSE}_{\Delta x}/{\rm RMSE}_{\Delta x/2}\big)}{\log 2}.
	\end{align}
	
	In addition, to preserve the fourth-order accuracy of  the present MRT-LB model, the initial condition of distribution function $f_i(\bm{x},t)$ in the MRT-LB model must be given properly. Based on the previous work \cite{Chai2020}, we utilize the first two terms of Eq. (\ref{f-feq-o5}) to initialize the distribution function $f_i$ in the implementation of the MRT-LB model.\\
	\textbf{Example 1.} We first consider the one-dimensional Burgers' equation (\ref{eq-u}) with the initial and boundary conditions,
	\begin{align}\label{Ex-1D-Burgers-E}
		\begin{cases}
			u(x,0)=\sin(\pi x),&0\leq x\leq 1,\\
			u(0,t)=u(1,t)=0, &t>0,
		\end{cases}
	\end{align} 
	and the analytical solution of $u(x,t)$ is given by
	\begin{align}\label{soul-Ex-1D-Burgers-E}
		u(x,t)=4\pi\upsilon \frac{\sum_{n=1}^{+\infty}a_n\exp{\big(-n^2\pi^2\upsilon t\big)}n\sin(n\pi x)}{a_0+2\sum_{n=1}^{+\infty}a_n\exp{\big(-n^2\pi^2\upsilon t\big)}\cos(n\pi x)},
	\end{align}
	where $a_n=I_n\big(1/(2\pi \upsilon)\big)$ with
	%, $n=0,1,2,\ldots$, in which 
	$I_n(x)$ representing the first type of the $n$-th modified Bessel function. With the help of the Cole-Hopf transformation, one can derive the analytical solution of $\theta(x,t)$ in the one-dimensional diffusion equation,
	\begin{align}\label{soul-theta-Ex-1D-Burgers-E}
		\theta(x,t)=a_0+2\sum_{n=1}^{+\infty}a_n\exp{\big(-n^2\pi^2\upsilon t\big)}\cos(n\pi x).
	\end{align}
	%Here it should be noted that the upper bound of index $n$ in Eqs. (\ref{soul-Ex-1D-Burgers-E}) and (\ref{soul-theta-Ex-1D-Burgers-E}) should be larger enough in simulations, and  $N=10000$ is adopted.
	
\begin{figure}   % 常规操作\begin{figure}开头说明插入图片
			% 后面跟着的[htbp]是图片在文档中放置的位置，也称为浮动体的位置，关于这个我们后面的文章会聊聊，现在不管，照写就是了
			\centering            % 前面说过，图片放置在中间
			\subfloat[$\theta(x,t)$]    
			{
				\includegraphics[width=0.4\textwidth]{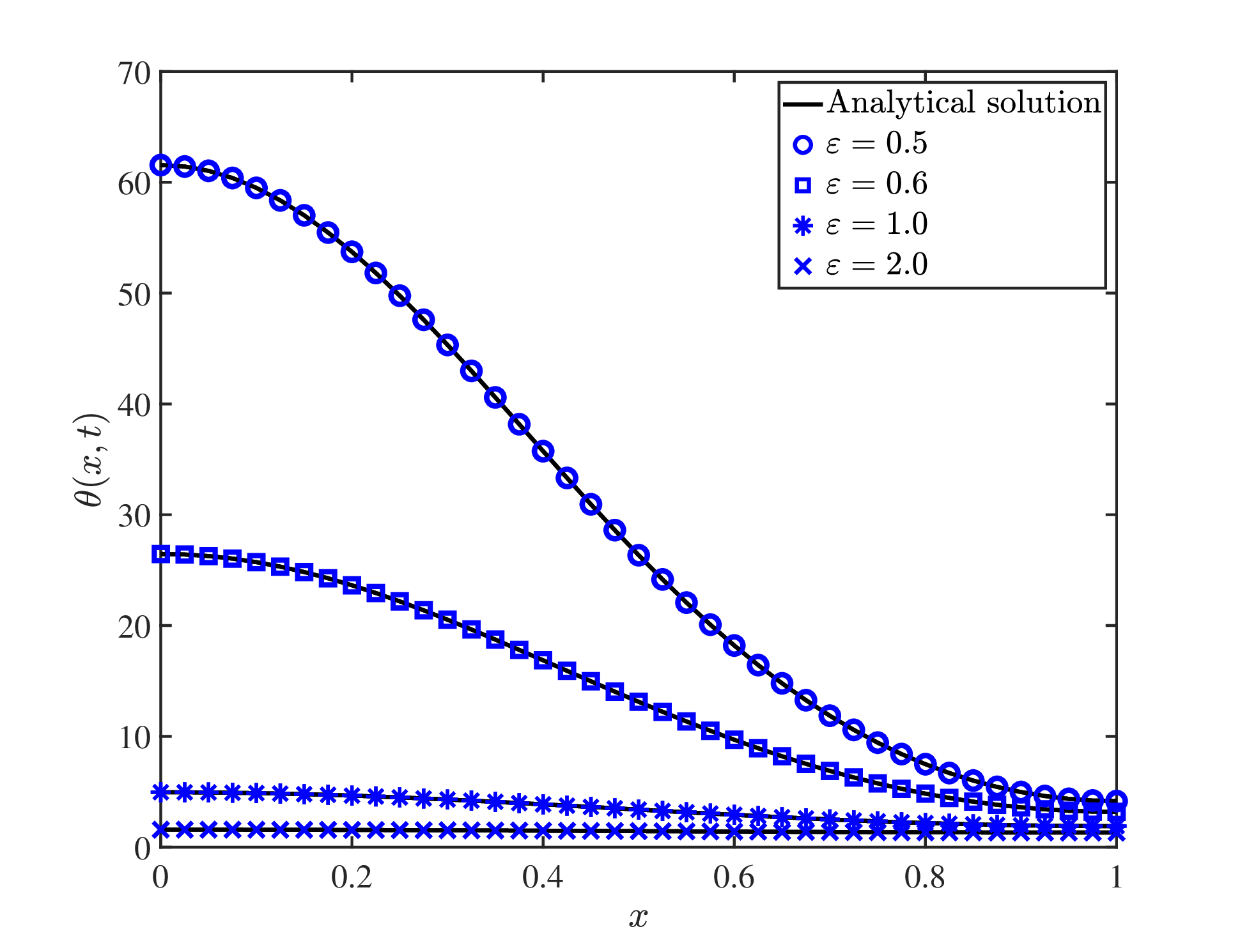} 
			}
			\subfloat[  $u(x,t)$]
			{
				\includegraphics[width=0.4\textwidth]{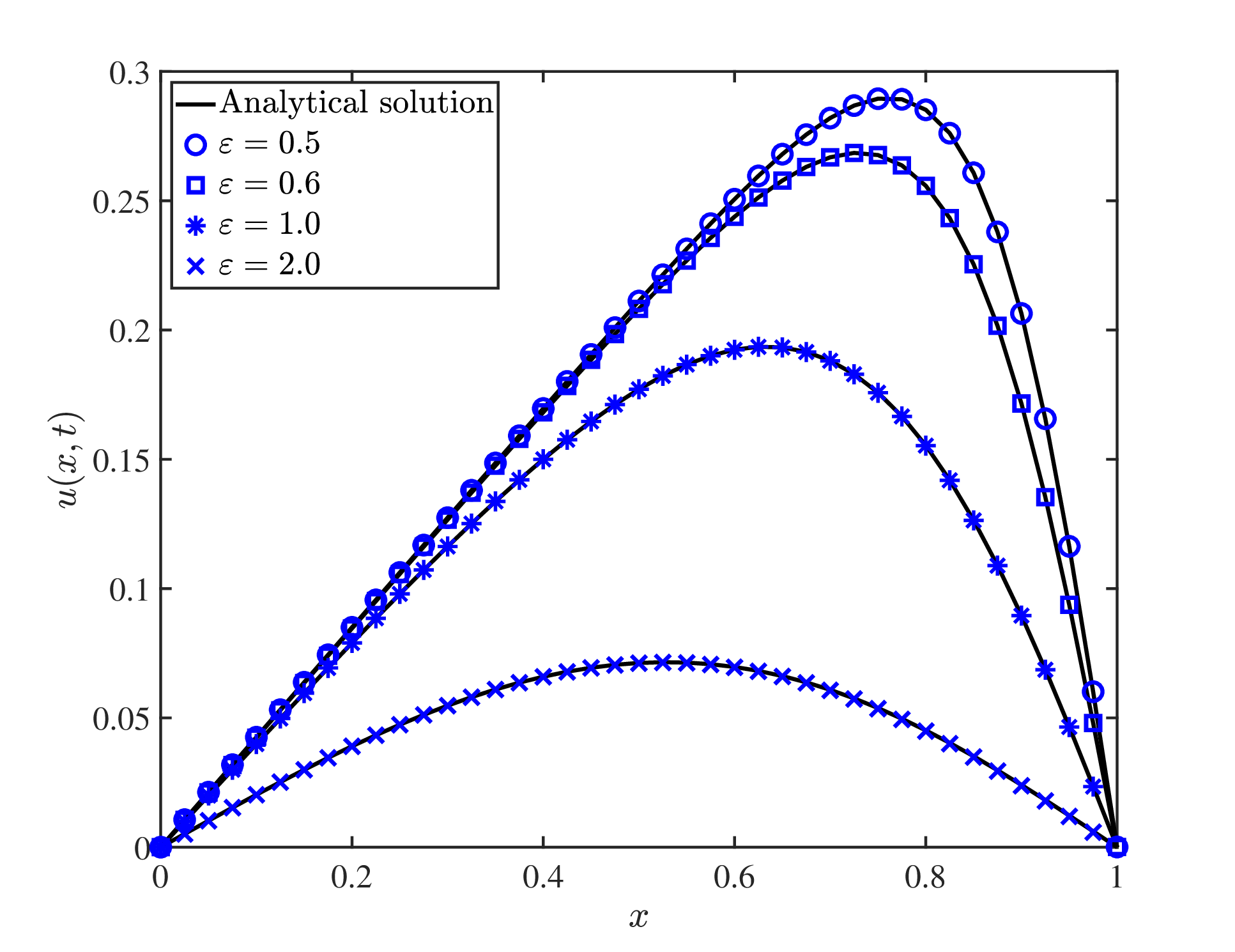}
			}
			\caption{Comparisons of the analytical and numerical results under different values of parameter  $\varepsilon$ ($t = 2.0$)  (solid
				lines: analytical solutions, symbols: numerical results).}    % 整个图片的说明，注释写在{}内
			\label{Fig-Ex1-diff-kappa}       % 整个图片的标签编号，注意这里跟子图是一样的道理，标签不能重复 
		\end{figure}

\begin{figure}    % 常规操作\begin{figure}开头说明插入图片
				% 后面跟着的[htbp]是图片在文档中放置的位置，也称为浮动体的位置，关于这个我们后面的文章会聊聊，现在不管，照写就是了
				\centering            % 前面说过，图片放置在中间
				\subfloat[$\theta(x,t)$]    
				{
					\includegraphics[width=0.4\textwidth]{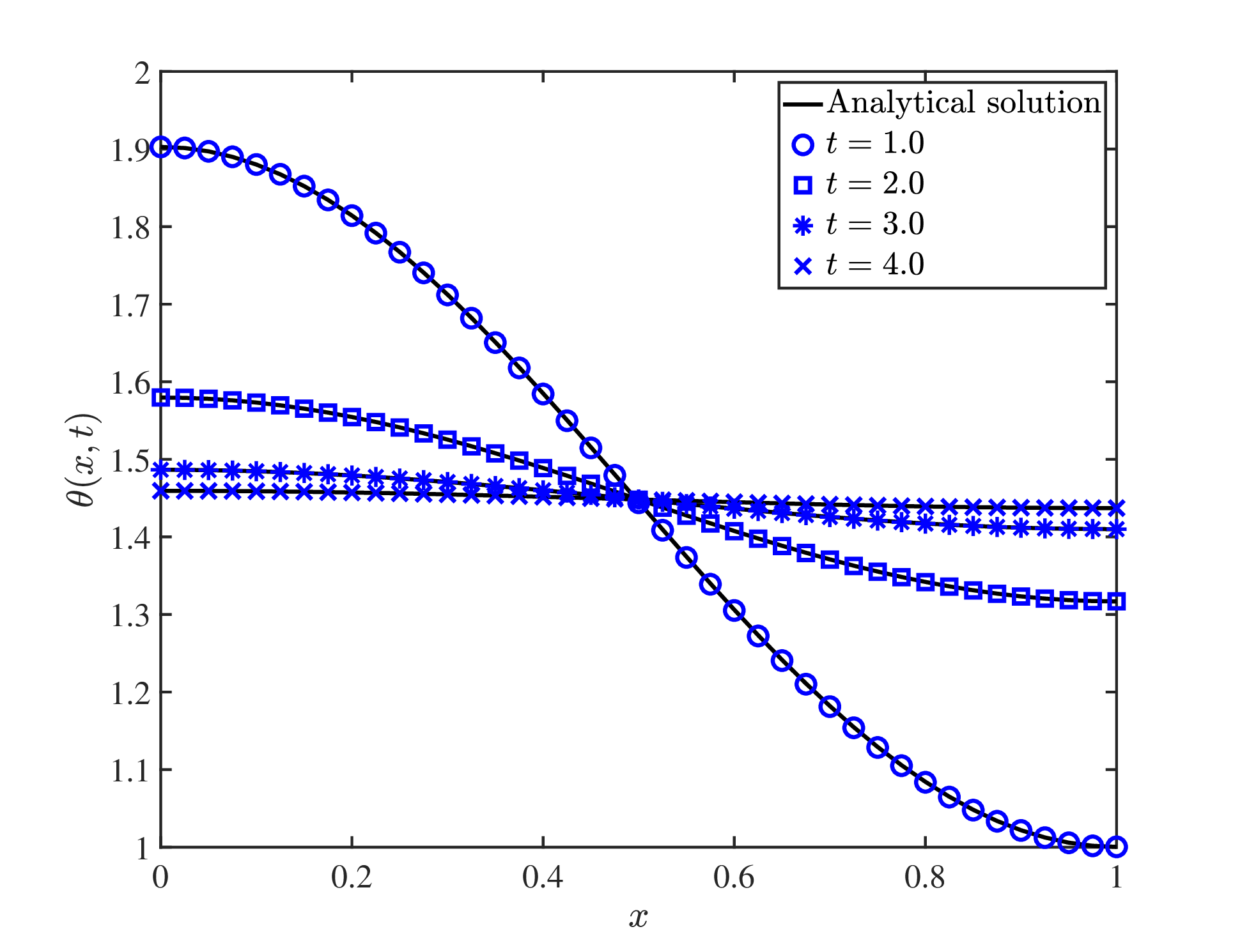} 
				}
				\subfloat[  $u(x,t)$]
				{
					\includegraphics[width=0.4\textwidth]{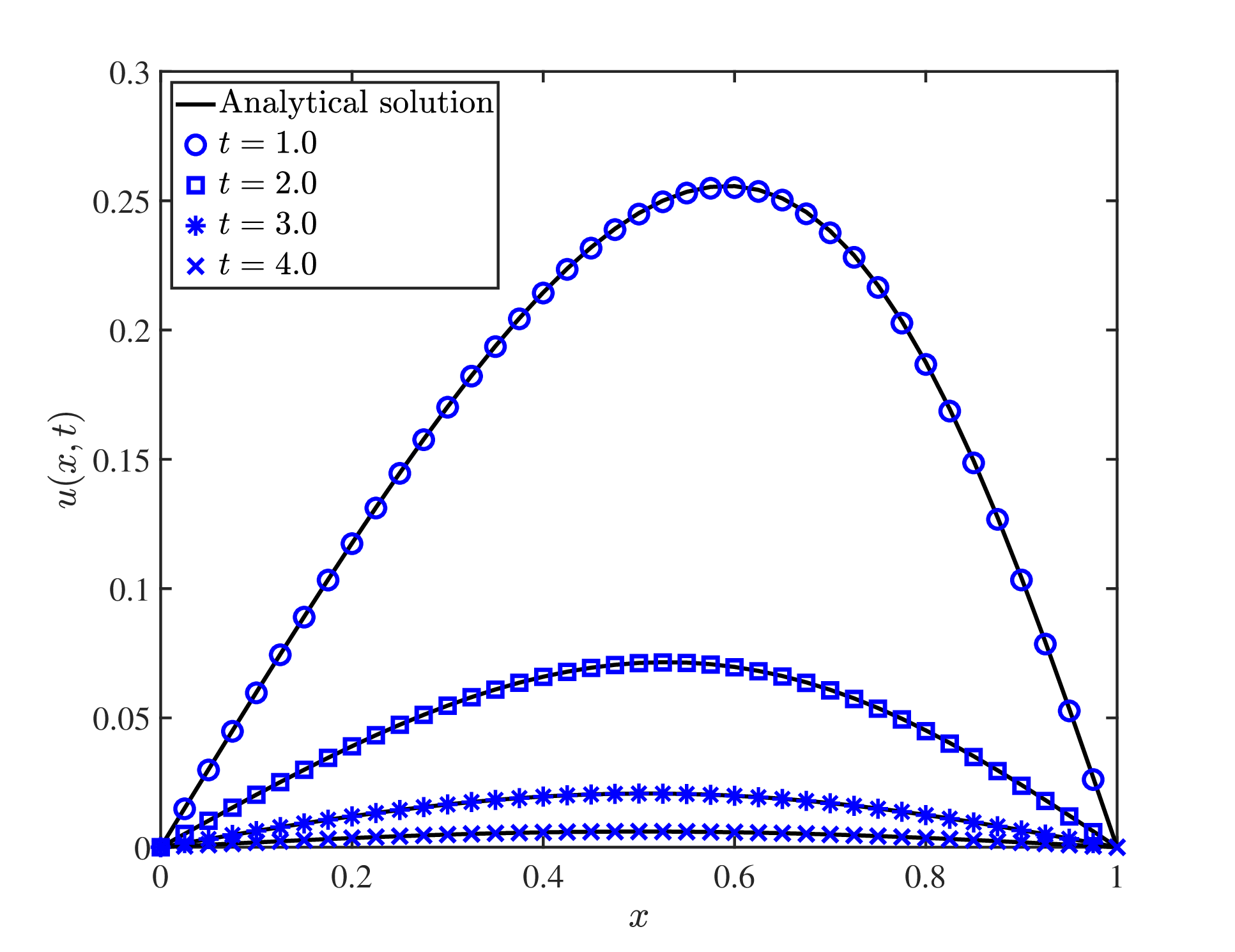}
				}
				\caption{Comparisons of the analytical and numerical results with parameter   at different values of   time   $t$ ($\varepsilon=2.0$)  (solid
					lines: analytical solutions, symbols: numerical results).}    % 整个图片的说明，注释写在{}内
				\label{Fig-Ex1-diff-t}       % 整个图片的标签编号，注意这里跟子图是一样的道理，标签不能重复 
			\end{figure}

	We first conduct some simulations under different values of parameter  $\varepsilon$ ($\varepsilon=$0.5, 0.6, 1.0, 2.0) corresponding to different viscosity coefficients for the specific lattice spacing $\Delta x=1/40$ and time step $\Delta t=1/100$, and plot the results in Fig. \ref{Fig-Ex1-diff-kappa}. Additionaly, to see the evolution of variables $\theta(x,t)$ and $u(x,t)$ in time, we also carry out some simulations at different values of time $t$ ($t =$ 1.0, 2.0, 3.0, 4.0), and present the results in Fig. \ref{Fig-Ex1-diff-t} where   $\varepsilon = 2.0$. As shown in
	Figs. \ref{Fig-Ex1-diff-kappa} and \ref{Fig-Ex1-diff-t},  the numerical results of MRT-LB model  are in good
	agreement with the corresponding analytical solutions.
			
			In order to measure the deviation between the numerical results and analytical solutions, we also perform some simulations under different values of lattice spacing $\Delta x$,  and list the RMSEs and   CRs of the MRT-LB model in Table \ref{Table-1D}. From the result shown in this table, one can observe that the MRT-LB model is fourth-order accurate in space, which is in agreement with our theoretical analysis.

			\begin{table} 
				\begin{center}
					\caption{The RMSEs and CRs under different values of  parameter  $\varepsilon$ ($t=2.0$).} 
					\label{Table-1D}
					\begin{tabular}{cccccccccccccc}\hline\hline
						$\varepsilon$&$\Delta x$&$\Delta t$&variable&RMSE$_{\Delta x,\Delta t}$&RMSE$_{\Delta x/2,\Delta t/4}$&RMSE$_{\Delta x/4,\Delta t/16}$&RMSE$_{\Delta x/8,\Delta t/64}$&CR\\
						\hline
						\multirow{2}{*}{0.5}&\multirow{2}{*}{1/40}&\multirow{2}{*}{1/100}&$\theta$&7.8613$\times 10^{-4}$&4.8926$\times 10^{-5}$&3.0595$\times 10^{-6}$&1.9139$\times 10^{-7}$&$\sim$4.0013\\ 
						&&&$u$&1.9827$\times 10^{-5}$&1.1821$\times 10^{-6}$&7.2429$\times 10^{-8}$&4.4865$\times 10^{-9}$&$\sim$4.0365\\ 
						\cline{4-9}
						\multirow{2}{*}{1}&\multirow{2}{*}{1/40}&\multirow{2}{*}{1/100}&$\theta$& 2.4232$\times 10^{-4 }$& 1.5014$\times 10^{-5 }$&9.3770 $\times 10^{-7 }$& 5.8641$\times 10^{- 8}$&$\sim$4.0042\\ 
						&&&$u$&3.2184 $\times 10^{-5}$&2.0194$\times 10^{-6}$&1.2608$\times 10^{-7}$&7.8713$\times 10^{-9}$&$\sim$3.9992\\ 
						\cline{4-9}
						\multirow{2}{*}{1.5}&\multirow{2}{*}{1/40}&\multirow{2}{*}{1/100}&$\theta$&2.3106$\times 10^{-4}$&1.3918 $\times 10^{-5}$&8.6303$\times 10^{-7}$&5.3874$\times 10^{-8}$&$\sim$4.0221\\ 
						&&&$u$&1.0796$\times 10^{-4}$&6.4468$\times 10^{-6}$&3.9752$\times 10^{-7}$&2.4739$\times 10^{-8}$&$\sim$4.0305\\ 
						\cline{4-9}
						\multirow{2}{*}{2}&\multirow{2}{*}{1/40}&\multirow{2}{*}{1/100}&$\theta$&2.1330$\times 10^{-4}$&1.2349$\times 10^{-5}$&7.5838$\times 10^{-7}$&4.7229$\times 10^{-8}$&$\sim$4.0470\\ 
						&&&$u$&1.8574$\times 10^{-4}$&1.0638$\times 10^{-5}$&6.4967$\times 10^{-7}$&4.0343$\times 10^{-8}$&$\sim$4.0562\\ 
						\hline\hline
					\end{tabular}
				\end{center} 
			\end{table}
			\noindent\textbf{Example 2.} We continue to consider the two-dimensional coupled Burgers' equations (\ref{eq-u}) with the periodic boundary condition and the following initial condition,
			\begin{align}\label{Ex-2D-Burgers-E}
				\begin{cases}
					u_x(x,y,0)=-2\upsilon\frac{ 2\pi \cos(2\pi x)\sin(\pi y) }{2+\sin(2\pi x)\sin(\pi y)},&(x,y)\in\Omega,\\
					u_y(x,y,0)=-2\upsilon\frac{ \pi  \sin(2\pi x)\cos(\pi y)}{2+\sin(2\pi x)\sin(\pi y) } ,&(x,y)\in\Omega,
				\end{cases}
			\end{align} 
			where the computational domain is $\Omega=\{(x,y):0\leq x, y\leq 2\}$. Under the condition of Eq. (\ref{Ex-2D-Burgers-E}), one can obtain  the analytical solution  of $\bm{u}=(u_x,u_y)^T$,
			\begin{align}\label{soul-Ex-2D-Burgers-E}
				\left\{\begin{aligned}
					&u_x(x,y,t)=-2\upsilon\frac{ 2\pi \exp{\big(-5\upsilon \pi^2 t\big)}\cos(2\pi x)\sin(\pi y) }{2+\exp{\big(-5\upsilon \pi^2 t\big)}\sin(2\pi x)\sin(\pi y)},\\
					&u_y(x,y,t)=-2\upsilon\frac{ \pi \exp{\big(-5\upsilon \pi^2 t\big)}\sin(2\pi x)\cos(\pi y)}{2+\exp{\big(-5\upsilon \pi^2 t\big)}\sin(2\pi x)\sin(\pi y) },
				\end{aligned}\right.
			\end{align} 
			then according to Eq. (\ref{eq-theta}), the analytical solution  of  $\theta(x,y,t)$ can be derived   \cite{Zhanlav2016},
			\begin{align}\label{soul-theta-Ex-2D-Burgers-E}
				\theta(x,y,t)=\frac{2+\exp{\big(-5\upsilon \pi^2 t\big)}\sin(2\pi x)\sin(\pi y)}{2}.
			\end{align}
			
\begin{figure}    % 常规操作\begin{figure}开头说明插入图片
					% 后面跟着的[htbp]是图片在文档中放置的位置，也称为浮动体的位置，关于这个我们后面的文章会聊聊，现在不管，照写就是了
					\centering            % 前面说过，图片放置在中间
					\subfloat[$\theta(x,y,t)$]    
					{
						\includegraphics[width=0.4\textwidth]{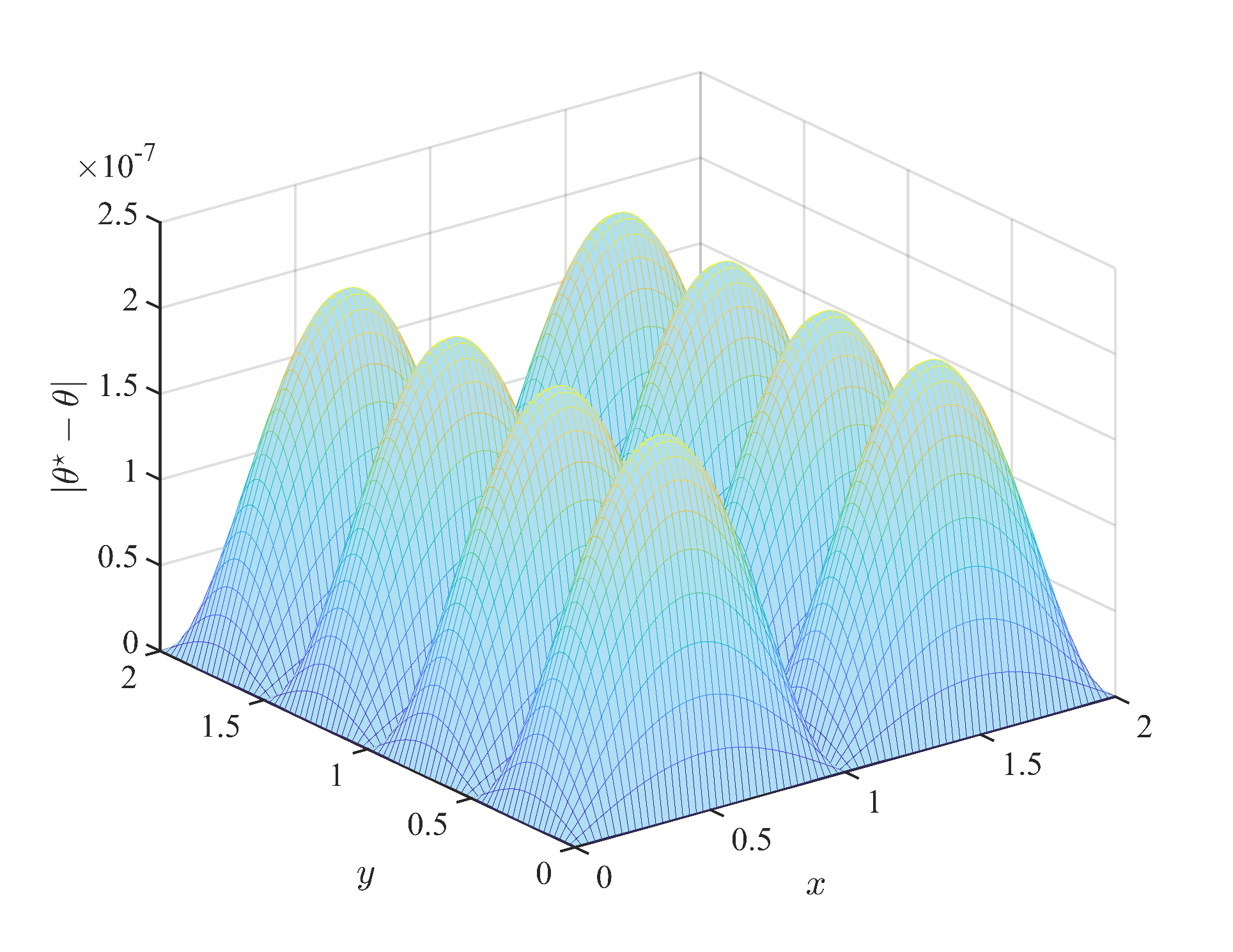} 
					}
					\subfloat[  $u_x(x,y,t)$]
					{
						\includegraphics[width=0.4\textwidth]{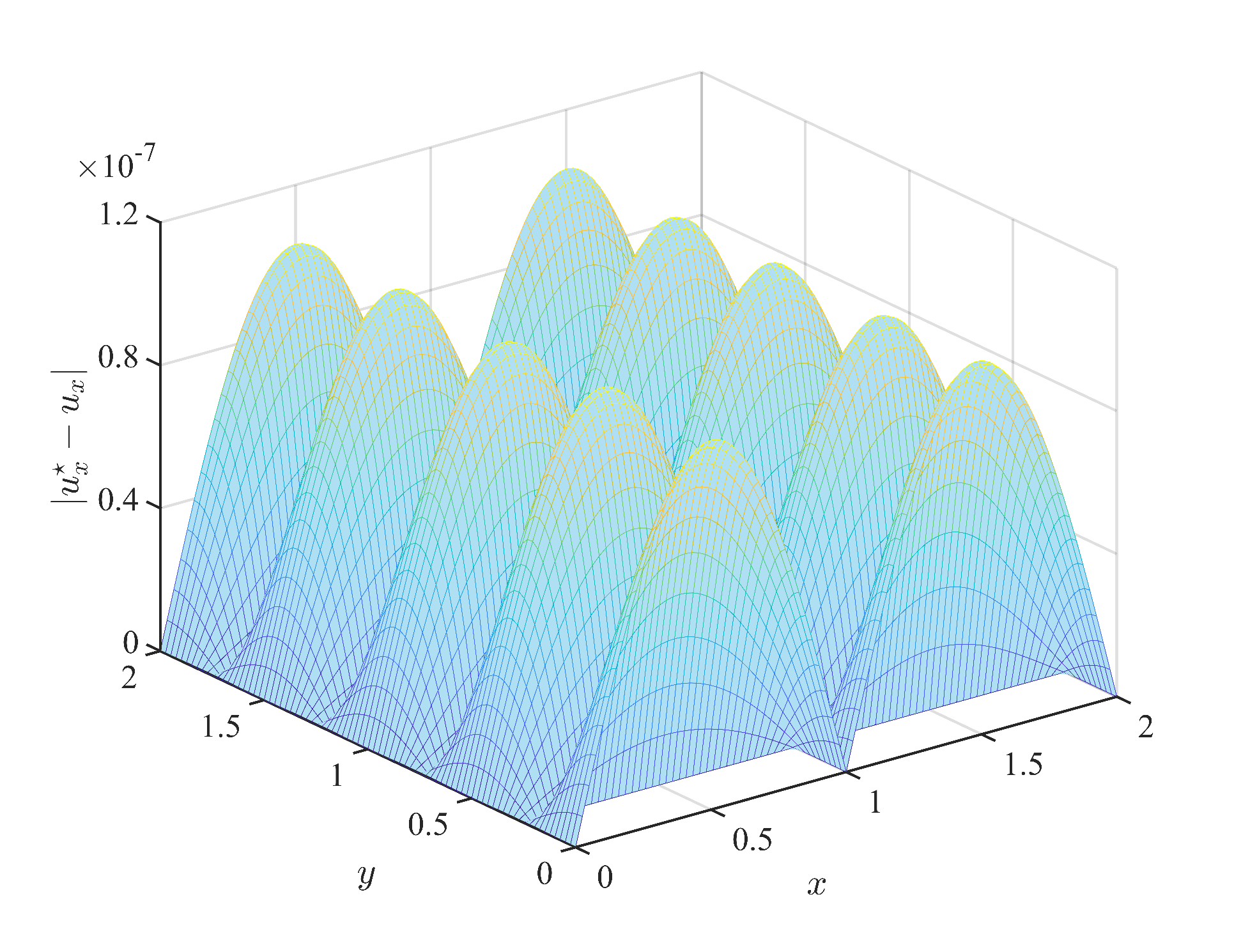}
					}
					\quad
					
					\subfloat[  $u_y(x,y,t)$]
					{
						\includegraphics[width=0.4\textwidth]{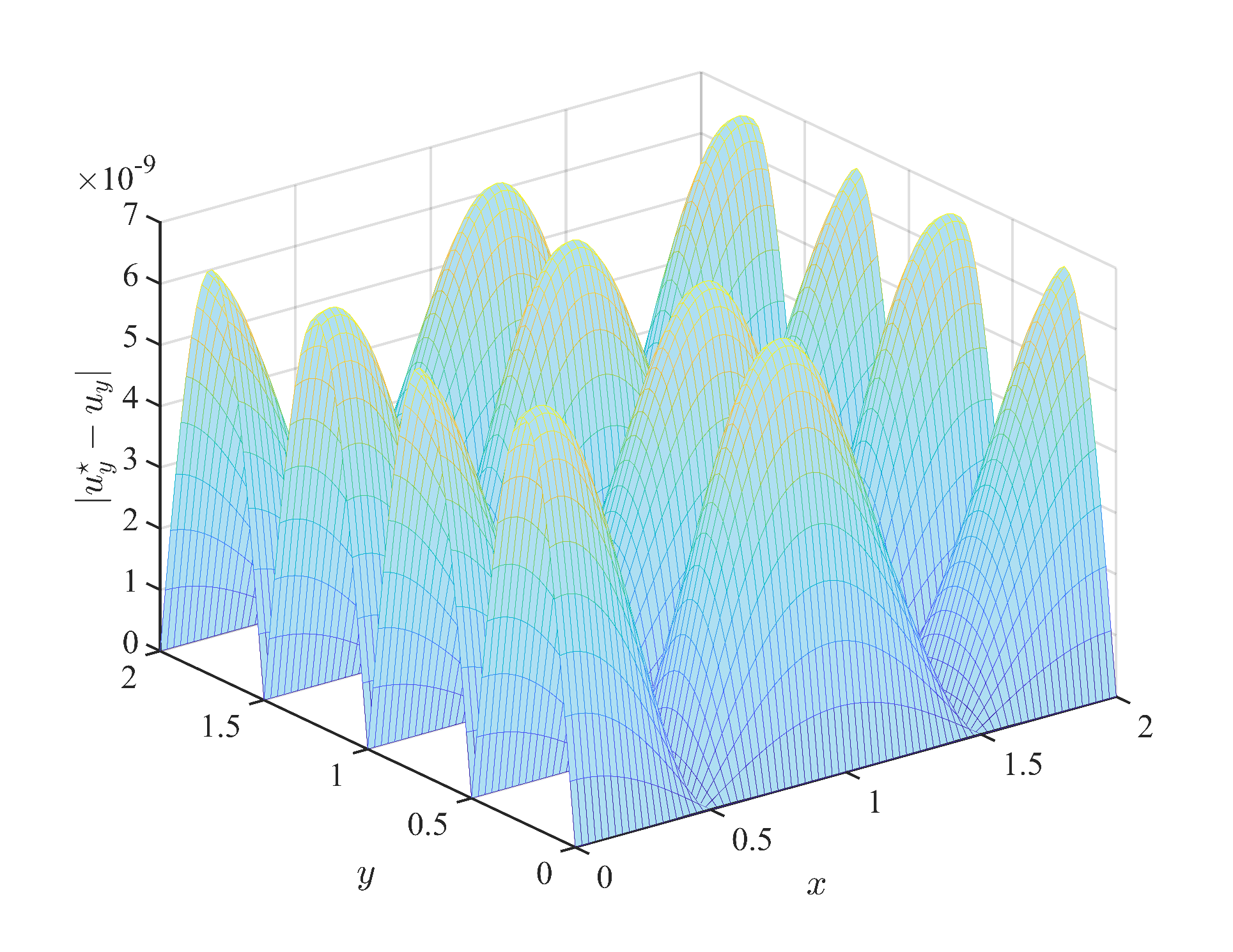}
					}
					\caption{The absolute errors between the analytical and numerical results of variables $\theta$, $u_x$  and $u_y$ ($\varepsilon=0.2$, $t=2.0$).}    % 整个图片的说明，注释写在{}内
					\label{Fig-Ex-2D-error}     % 整个图片的标签编号，注意这里跟子图是一样的道理，标签不能重复 
				\end{figure}

\begin{figure}    % 常规操作\begin{figure}开头说明插入图片
						% 后面跟着的[htbp]是图片在文档中放置的位置，也称为浮动体的位置，关于这个我们后面的文章会聊聊，现在不管，照写就是了
						\centering            % 前面说过，图片放置在中间
						\subfloat[$\theta(x,1/4,t)$]    
						{
							\includegraphics[width=0.4\textwidth]{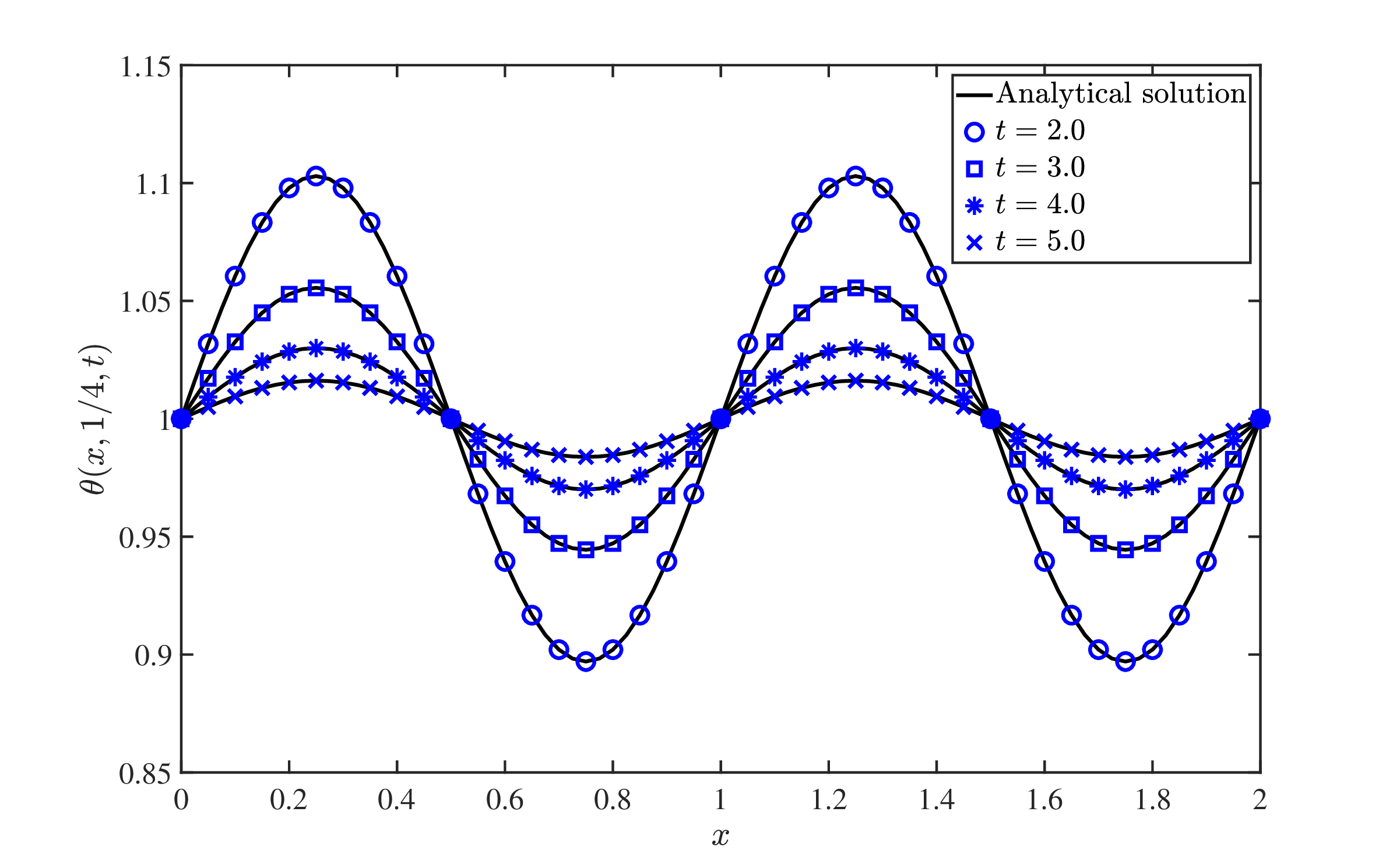} 
						}
						\subfloat[  $u_x(x,1/4,t)$]
						{
							\includegraphics[width=0.4\textwidth]{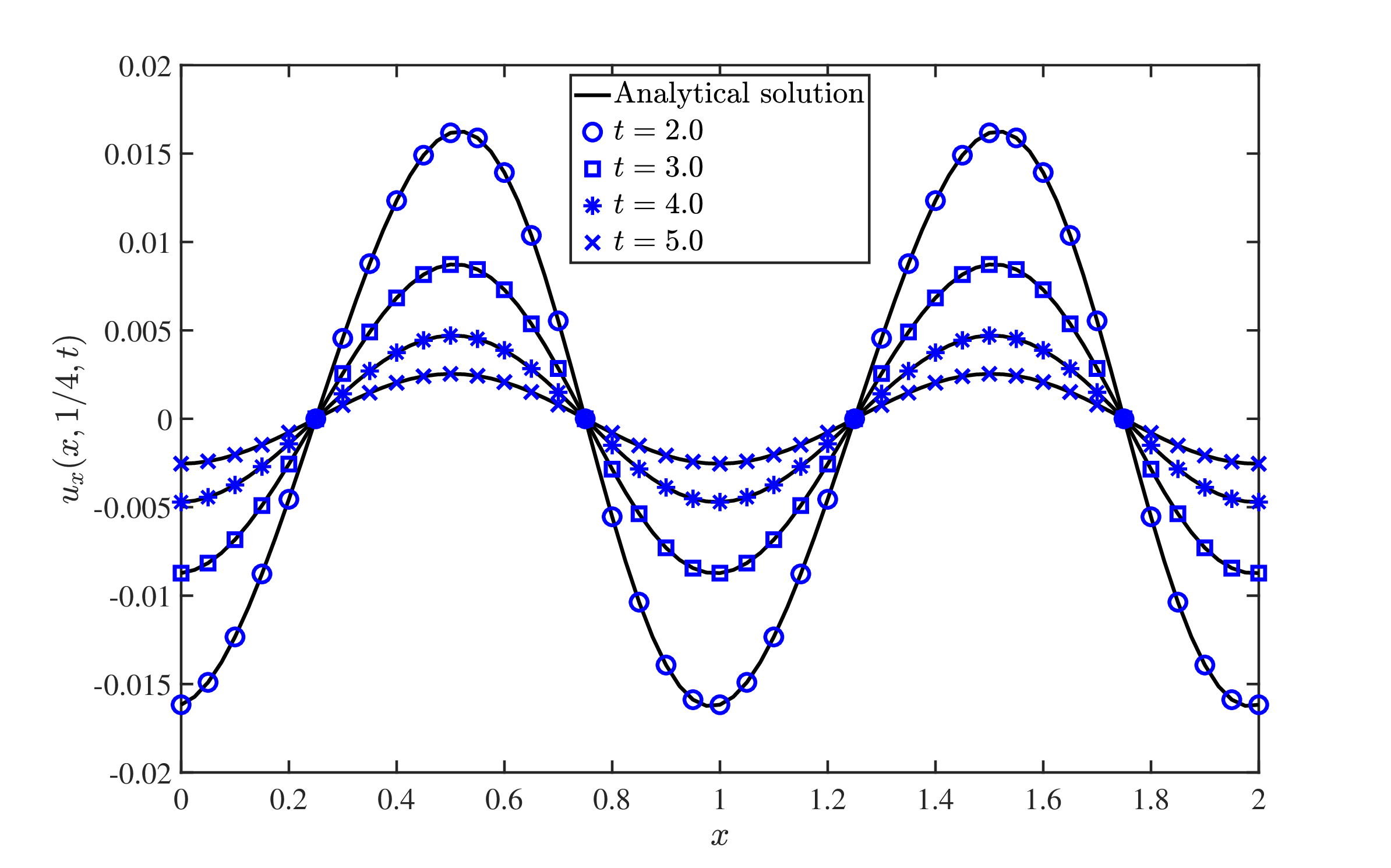}
						}
						\quad
						
						\subfloat[  $u_y(x,1/4,t)$]
						{
							\includegraphics[width=0.4\textwidth]{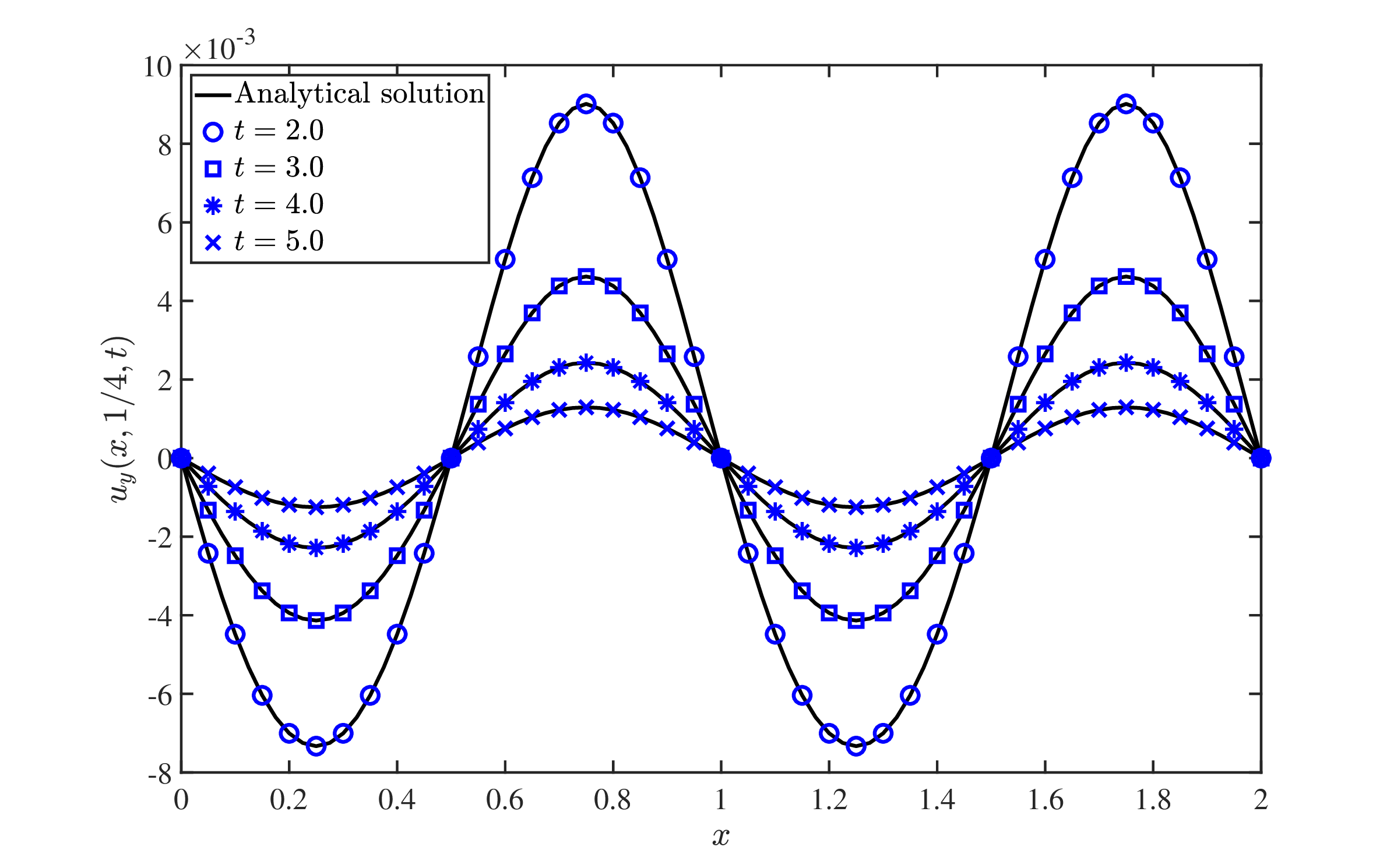}
						}
						\caption{Profiles of the variables at different values of time $t$ ($\varepsilon=0.2$) (solid
							lines: analytical solutions, symbols: numerical results).}    % 整个图片的说明，注释写在{}内
						\label{Fig-Ex-2D-pou}     % 整个图片的标签编号，注意这里跟子图是一样的道理，标签不能重复 
					\end{figure}

			First of all, for the given lattice spacing $\Delta x=1/40$ and time step $\Delta t=1/100$, we set the parameter $\varepsilon=0.2$ to determine the  viscosity coefficient $\upsilon$, and measure the absolute errors between the analytical and numerical results ($|\psi^{\star}-\psi|$, $\psi=\theta,u_x,u_y$) of the variables  $\theta$, $u_x$ and $u_y$   in  Fig. \ref{Fig-Ex-2D-error}  where $t=2.0$. As shown in this figure,  the  maximum absolute errors of the variables $\theta$, $u_x$  and $u_y$ are less than $2.5\times 10^{-7}$, $1.2\times 10^{-7}$  and $7.0\times 10^{-9}$, respectively. Then we also conduct some simulations at
			different values of time $t$ ($t=$2.0, 3.0, 4.0, 5.0), and plot the profiles of the variables $\theta$, $u_x$  and $u_y$ in Fig.  \ref{Fig-Ex-2D-pou} where $\varepsilon=0.2$. As seen from this figure, the numerical results are very close to the analytical solutions. Finally, to test the CR  of the MRT-LB model, we measure the RMSEs between the numerical and analytical solutions, and calculate the average CRs under  different values of  parameter $\varepsilon$ in Table \ref{Table-2D} where  $t=2.0$. From this table, one can observe that the MRT-LB model for two-dimensional coupled Burgers' equations also have a fourth-order convergence rate in space.
					
					\begin{table}
						\begin{center}
							\caption{The RMSEs and CRs under different values of  parameter  $\varepsilon$ ($t=2.0$).} 
							\label{Table-2D}
							\begin{tabular}{cccccccccccccc}\hline\hline
								$\varepsilon$&$\Delta x$&$\Delta t$&variable&RMSE$_{\Delta x,\Delta t}$&RMSE$_{\Delta x/2,\Delta t/4}$&RMSE$_{\Delta x/4,\Delta t/16}$&RMSE$_{\Delta x/8,\Delta t/64}$&CR\\
								\hline
								\multirow{3}{*}{0.1}&\multirow{3}{*}{1/20}&\multirow{3}{*}{1/50}&$\theta$&4.2630$\times 10^{-7}$&2.6622$\times 10^{-8}$&1.6687$\times 10^{-9}$&1.0453$\times 10^{-10}$&$\sim$3.9979    \\ 
								&&&$u_x$&2.6683$\times 10^{-7}$&1.6587$\times 10^{-8}$&1.0353$\times 10^{-9}$&6.4689$\times 10^{-11}$&$\sim$4.0034   \\ 
								&&&$u_y$&3.8117$\times 10^{-8}$&2.3404$\times 10^{-9}$&1.4282$\times 10^{-10}$&8.9076$\times 10^{-12}$&$\sim$ 4.0210 \\ 
								\cline{4-9}
								\multirow{3}{*}{0.2}&\multirow{3}{*}{1/20}&\multirow{3}{*}{1/50}&$\theta$& 1.5114$\times 10^{-6}$& 9.4845$\times 10^{-8}$&5.9523$\times 10^{-9}$&3.7300$\times 10^{-10}$&$\sim$  3.9948    \\
								&&&$u_x$&	1.5218$\times 10^{-7}$&9.2774$\times 10^{-9}$&5.7625$\times 10^{-10}$&3.5958$\times 10^{-11}$&$\sim$ 4.0157 \\ 
								&&&$u_y$&		1.3577$\times 10^{-7}$&8.4954$\times 10^{-9}$&5.3110$\times 10^{-10}$&3.3196$\times 10^{-11}$&	$\sim$3.9993 \\ 
								\cline{4-9}
								\multirow{3}{*}{0.3}&\multirow{3}{*}{1/20}&\multirow{3}{*}{1/50}&$\theta$&  5.7600$\times 10^{-6}$&3.5539$\times 10^{-7}$&2.2209$\times 10^{-8}$&1.3902$\times 10^{-9}$&$\sim$4.0055   \\ 
								&&&$u_x$&
								1.3681$\times 10^{-6}$&8.4242$\times 10^{-8}$&5.2444$\times 10^{-9}$&3.2745$\times 10^{-10}$&$\sim$ 4.0095  \\
								&&&$u_y$& 
								7.7594$\times 10^{-7}$&4.8271$\times 10^{-8}$&3.0124$\times 10^{-9}$&1.8820$\times 10^{-10}$&$\sim$4.0032\\
								\cline{4-9}
								\multirow{3}{*}{0.4}&\multirow{3}{*}{1/20}&\multirow{3}{*}{1/50}&$\theta$&    9.9386$\times 10^{-6}$&5.9304$\times 10^{-7}$& 3.6764$\times 10^{-8}$&2.2967$\times 10^{-9}$&$\sim$ 4.0264   \\ 
								&&&$u_x$& 4.4801$\times 10^{-6}$&2.6578$\times 10^{-7}$&1.6396$\times 10^{-8}$&1.0214$\times 10^{-9}$&$\sim$ 4.0329  	\\
								&&&$u_y$& 
								2.1821$\times 10^{-6}$&1.3113$\times 10^{-7}$&8.1137$\times 10^{-9}$&5.0581$\times 10^{-10}$&$\sim$4.0249 \\ 
								\hline\hline
							\end{tabular}
						\end{center} 
					\end{table} 
					\noindent \textbf{Example 3.} We now focus on the three-dimensional coupled Burgers' equations (\ref{eq-u}) with the periodic condition and the following initial conditions,
					\begin{align}\label{Ex-3D-Burgers-E}
						\begin{cases}
							u_x(x,y,z,0)=-2\upsilon\frac{ 2\pi \cos(2\pi x)\sin(\pi y)\sin(4\pi z)}{2+\sin(2\pi x)\sin(\pi y)\sin(4\pi z)},&(x,y,z)\in\Omega,\\
							u_y(x,y,z,0)=-2\upsilon\frac{ \pi \sin(2\pi x)\cos(\pi y)\sin(4\pi z)}{2+ \sin(2\pi x)\sin(\pi y)\sin(4\pi z)},&(x,y,z)\in\Omega,\\
							u_z(x,y,z,0)=-2\upsilon\frac{ 4\pi  \sin(2\pi x)\sin(\pi y)\cos(4\pi z)}{2+ \sin(2\pi x)\sin(\pi y)\sin(4\pi z)},&(x,y,z)\in\Omega,
						\end{cases}
					\end{align} 
					where the computational domain is $\Omega=\{(x,y,z):0\leq x,y,z\leq 2\}$. The analytical solution of $\bm{u}=(u_x,u_y,u_z)^T$ is given by
					\begin{align}\label{soul-Ex-3D-Burgers-E}
						\left\{\begin{aligned}
							&u_x=-2\upsilon\frac{ 2\pi \exp{\big(-21\upsilon \pi^2 t\big)}\cos(2\pi x)\sin(\pi y)\sin(4\pi z)}{2+\exp{\big(-21\upsilon \pi^2 t\big)}\sin(2\pi x)\sin(\pi y)\sin(4\pi z)},\\
							&u_y=-2\upsilon\frac{ \pi \exp{\big(-21\upsilon \pi^2 t\big)}\sin(2\pi x)\cos(\pi y)\sin(4\pi z)}{2+\exp{\big(-21\upsilon \pi^2 t\big)}\sin(2\pi x)\sin(\pi y)\sin(4\pi z)},\\
							&u_z=-2\upsilon\frac{ 4\pi \exp{\big(-21\upsilon \pi^2 t\big)}\sin(2\pi x)\sin(\pi y)\cos(4\pi z)}{2+\exp{\big(-21\upsilon \pi^2 t\big)}\sin(2\pi x)\sin(\pi y)\sin(4\pi z)},
						\end{aligned}\right.
					\end{align} 
					from which one can also obtain the analytical solution of $\theta(x,y,z,t)$,
					\begin{align}
						&\theta(x,y,t)=\frac{2+\exp{\big(-21\upsilon \pi^2 t\big)}\sin(2\pi x)\sin(\pi y)\sin(4\pi z)}{2}.  
					\end{align}
					
\begin{figure}    % 常规操作\begin{figure}开头说明插入图片
							% 后面跟着的[htbp]是图片在文档中放置的位置，也称为浮动体的位置，关于这个我们后面的文章会聊聊，现在不管，照写就是了
							\centering            % 前面说过，图片放置在中间
							\subfloat[$\theta(x,y,1,t)$]    
							{
								\includegraphics[width=0.4\textwidth]{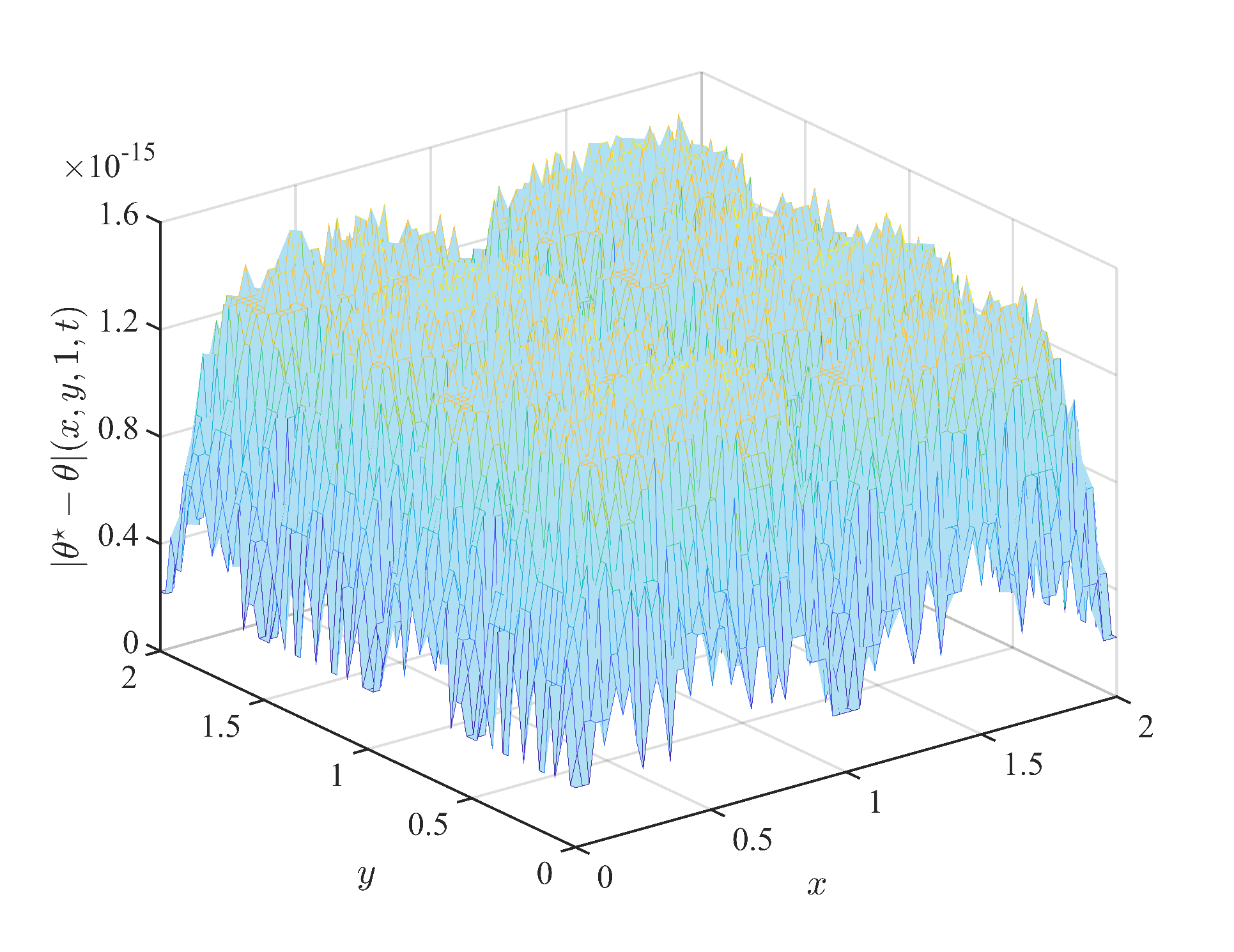} 
							}
							\subfloat[  $u_x(x,y,1,t)$]
							{
								\includegraphics[width=0.4\textwidth]{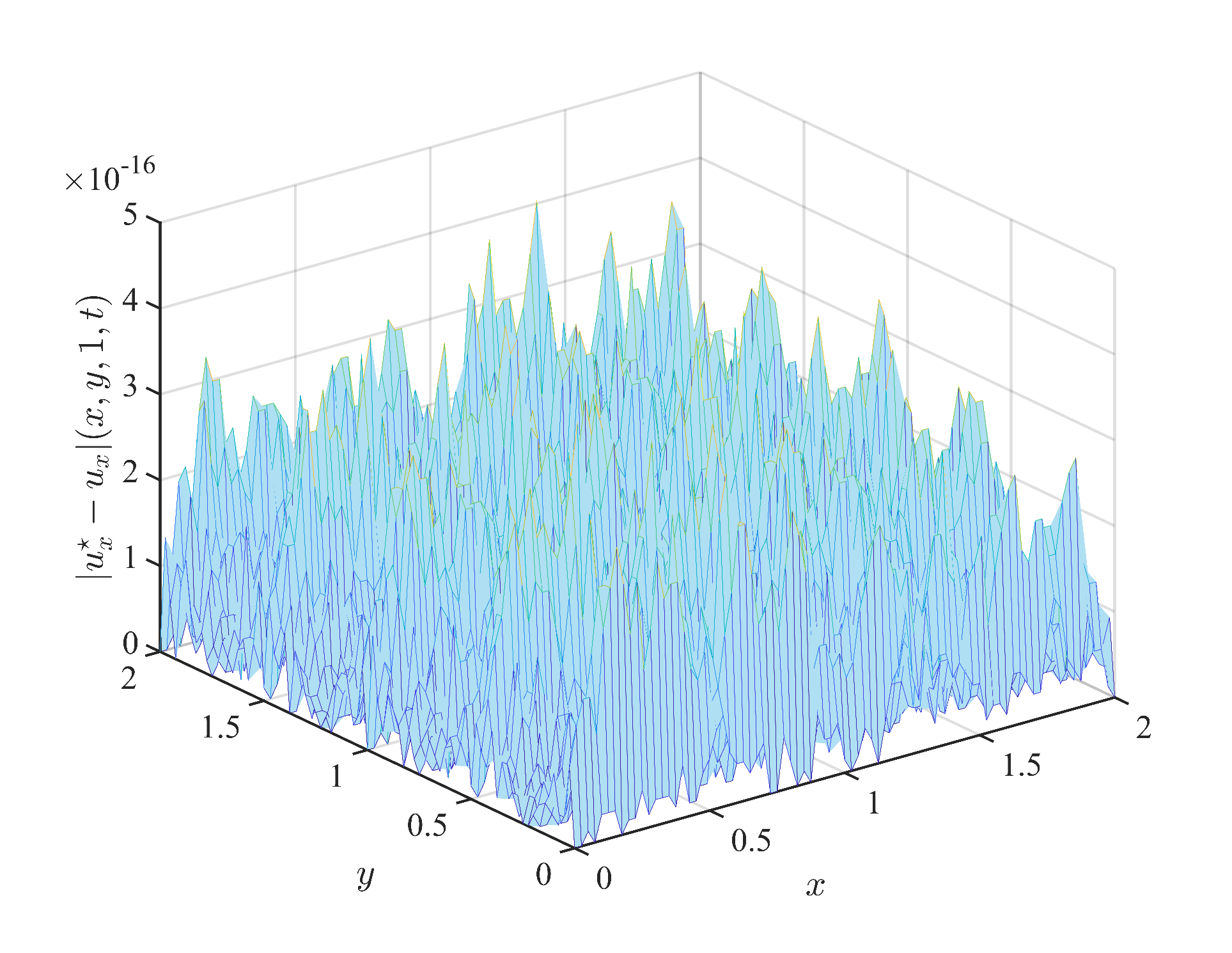}
							}
							\quad
							\subfloat[  $u_y(x,y,1,t)$]
							{
								\includegraphics[width=0.4\textwidth]{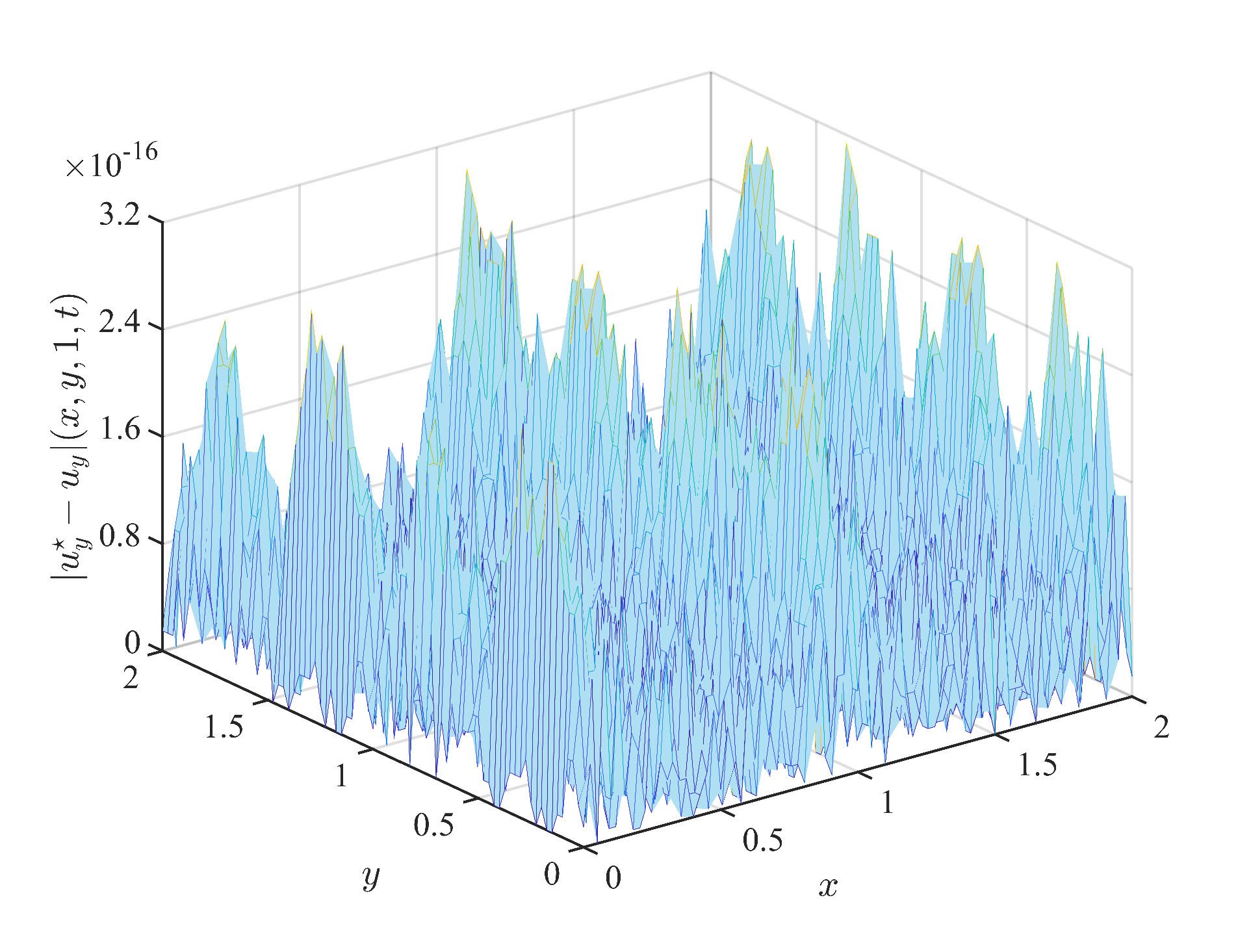}
							}
							\subfloat[  $u_z(x,y,1,t)$]
							{
								\includegraphics[width=0.4\textwidth]{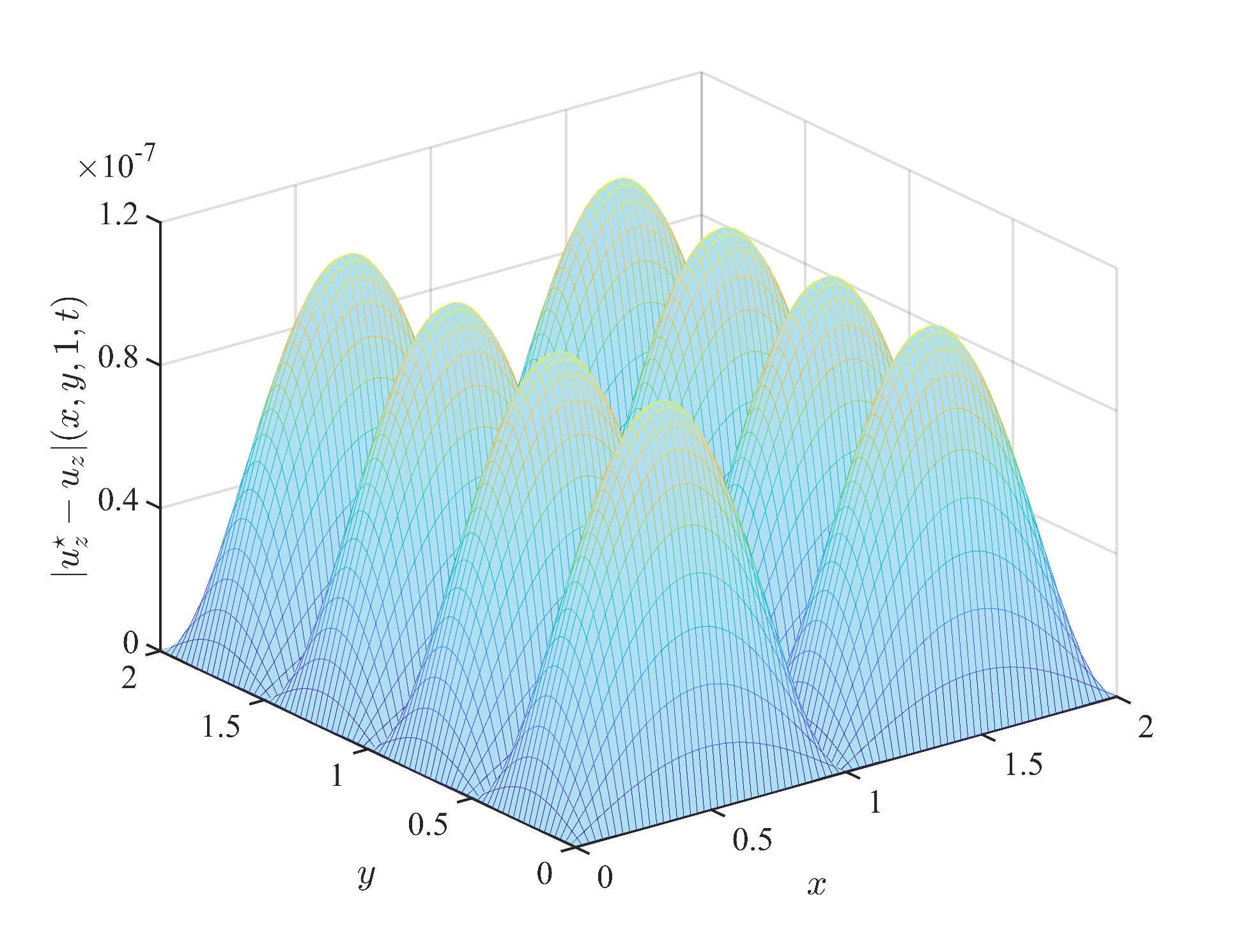}
							}
							\caption{The absolute errors between the analytical and numerical results of variables $\theta$, $u_x$, $u_y$  and $u_z$ ($\varepsilon=0.2$, $t=2.0$).}    % 整个图片的说明，注释写在{}内
							\label{Fig-Ex-3D-error}     % 整个图片的标签编号，注意这里跟子图是一样的道理，标签不能重复 
						\end{figure}

\begin{figure}    % 常规操作\begin{figure}开头说明插入图片
								% 后面跟着的[htbp]是图片在文档中放置的位置，也称为浮动体的位置，关于这个我们后面的文章会聊聊，现在不管，照写就是了
								\centering            % 前面说过，图片放置在中间
								\subfloat[$\theta(x,1/20,1/20,t)$]    
								{
									\includegraphics[width=0.4\textwidth]{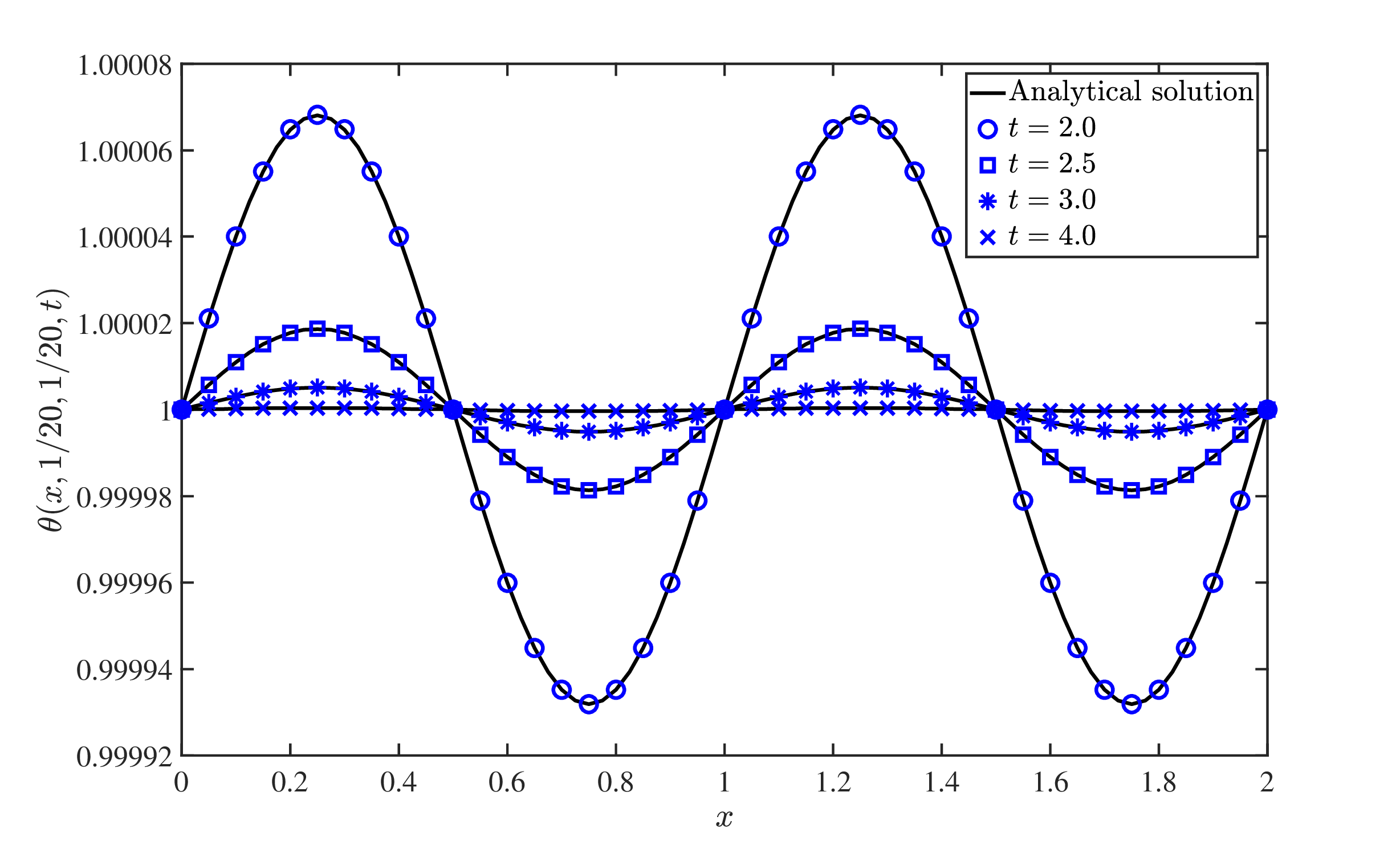} 
								}
								\subfloat[  $u_x(x,1/20,1/20,t)$]
								{
									\includegraphics[width=0.4\textwidth]{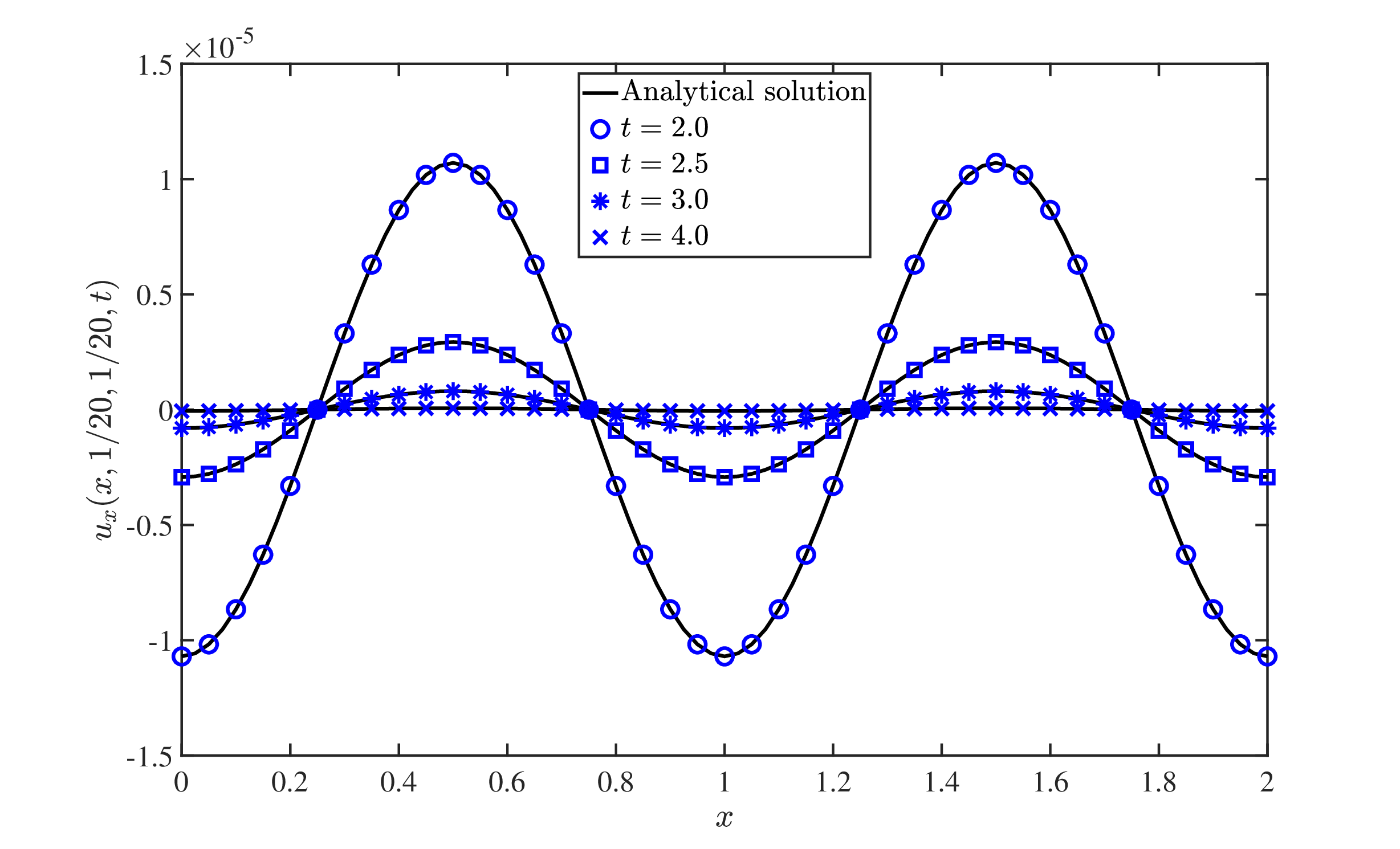}
								}
								\quad
								\subfloat[  $u_y(x,1/20,1/20,t)$]
								{
									\includegraphics[width=0.4\textwidth]{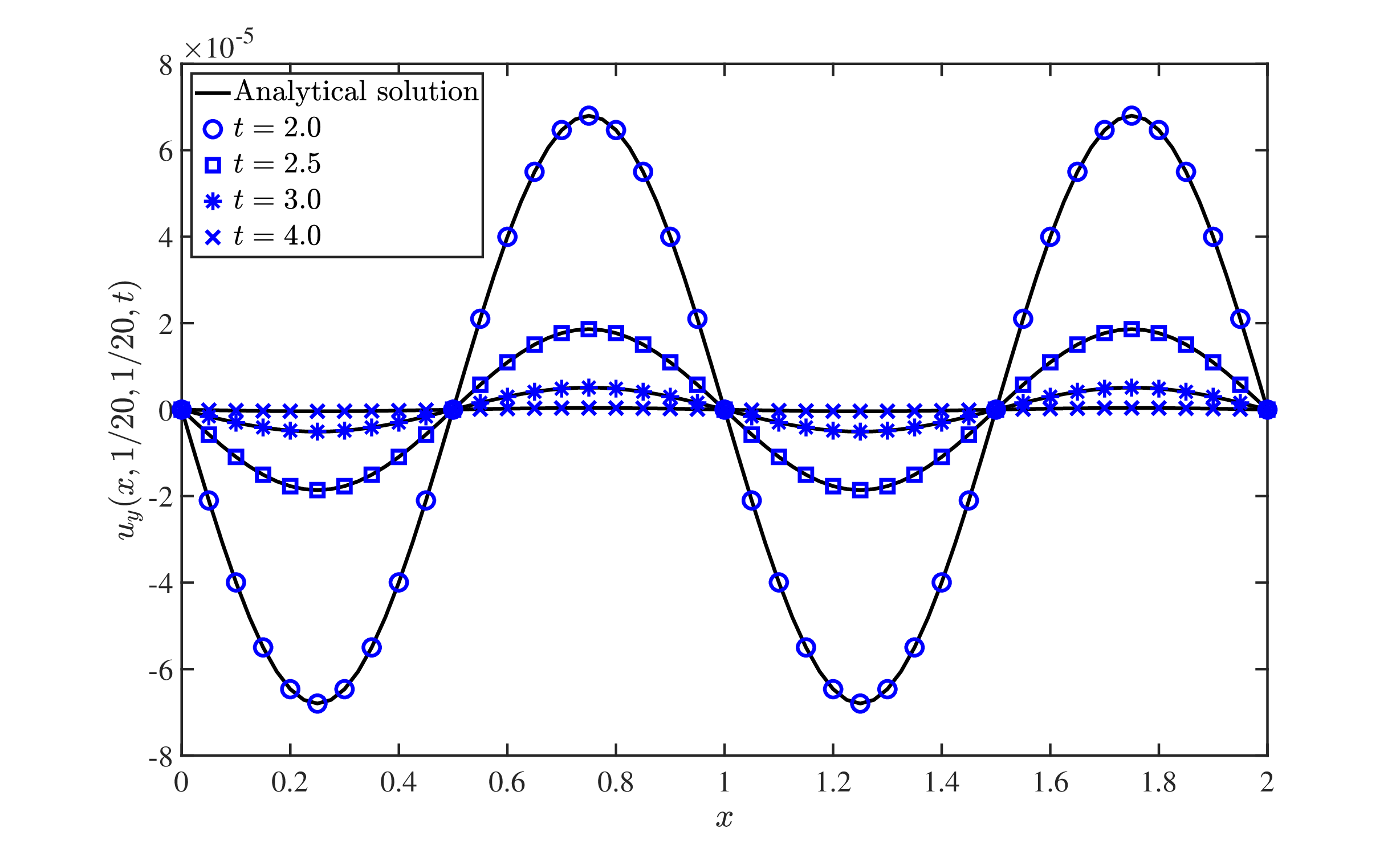}
								}
								\subfloat[  $u_z(x,1/20,1/20,t)$]
								{
									\includegraphics[width=0.4\textwidth]{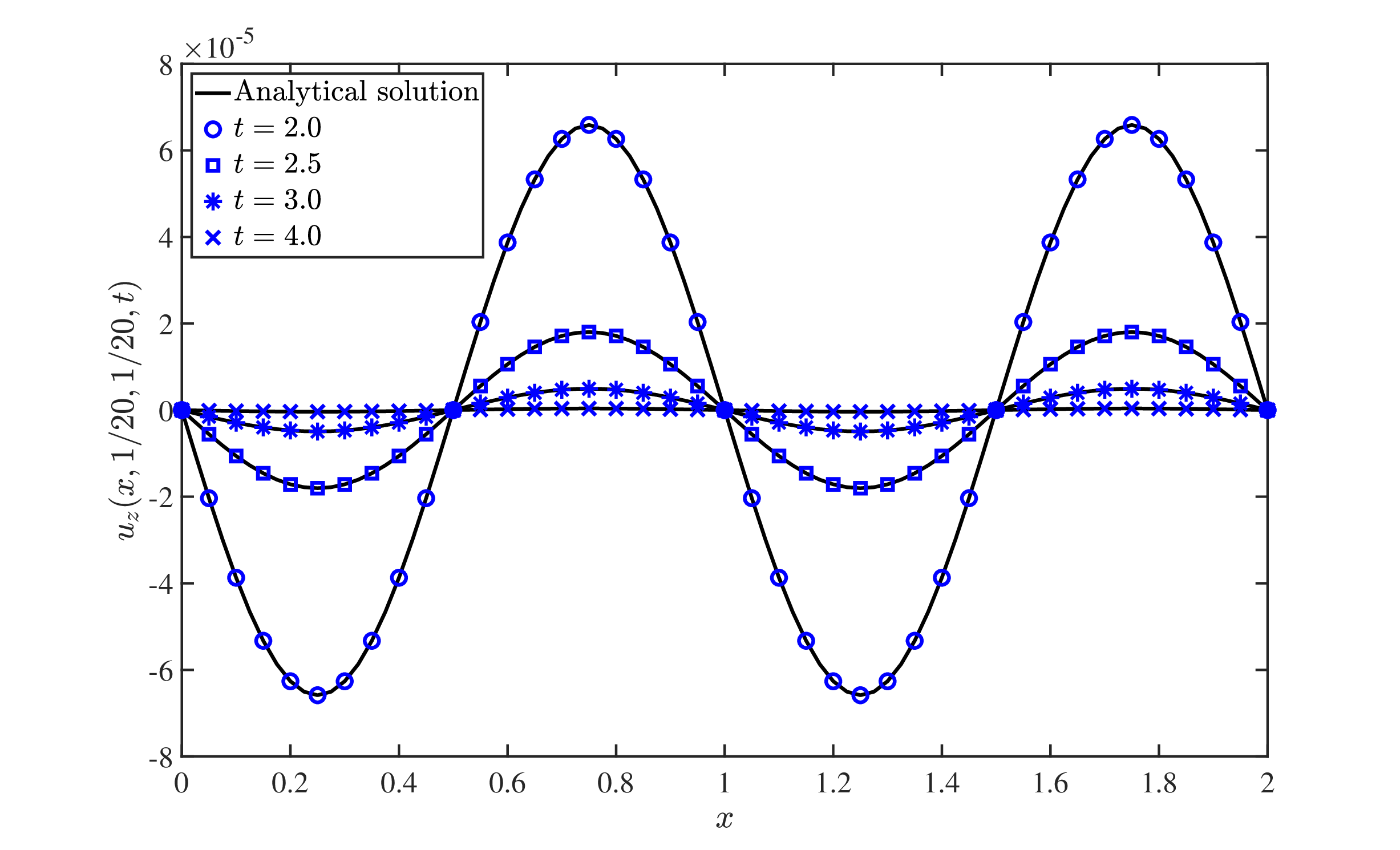}
								}
								\caption{Profiles of the variables at different values of time ($\varepsilon=0.2$) (solid
									lines: analytical solutions, symbols: numerical results).}    % 整个图片的说明，注释写在{}内
								\label{Fig-Ex-3D-pou}     % 整个图片的标签编号，注意这里跟子图是一样的道理，标签不能重复 
							\end{figure}

					Similar to the Example 2, in the following simulations, we fix the lattice spacing $\Delta x=1/40$ and time step $\Delta t=1/100$, and present the absolute errors between the analytical and numerical solutions of the variables $\theta$, $u_z$, $u_y$  and $u_z$ in  Fig. \ref{Fig-Ex-3D-error}  where   $\varepsilon=0.2$ and   $t=2.0$. From this figure, one can find that the maximum absolute errors of the variables $\theta$, $u_x$, $u_y$  and $u_z$ are less than $1.6\times 10^{-15}$, $5\times 10^{-15}$, $3.2\times 10^{-16}$  and $1.2\times 10^{-7}$, respectively. In order to see the evolution of variables $\theta$, $u_x$, $u_y$  and $u_z$  in time  more clearly, we further carry out some simulations at different values of time $t$ ($t$=2.0, 3.0, 4.0, 5.0) and $\varepsilon=0.2$, and plot the profiles of the variables $\theta$, $u_x$, $u_y$  and $u_z$   in Fig.  \ref{Fig-Ex-3D-pou}. As shown in this figure, the numerical results are in good agreement with the analytical solutions. Additionally, we also take different values of   parameter $\varepsilon$ ($\varepsilon=$0.05, 0.1, 0.15, 0.25) to calculate the RMSEs and CRs at the time $t=2.0$, and show the numerical results in Table \ref{Table-3D}  where a fourth-order accuracy in space is clearly observed.
					
							
							\begin{table} 
								\begin{center}
									\caption{The RMSEs and CRs under different values of parameter  $\varepsilon$  ($t=2.0$).} 
									\label{Table-3D}
									\begin{tabular}{cccccccccccccc}\hline\hline
										$\varepsilon$&$\Delta x$&$\Delta t$&variable&RMSE$_{\Delta x,\Delta t}$&RMSE$_{\Delta x/2,\Delta t/4}$&RMSE$_{\Delta x/4,\Delta t/16}$&RMSE$_{\Delta x/8,\Delta t/64}$&CR\\
										\hline
										\multirow{4}{*}{0.05}&\multirow{4}{*}{1/10}&\multirow{4}{*}{1/25}&$\theta$&1.4670$\times 10^{-5}$&7.9252$\times 10^{-7}$&4.7901$\times 10^{-8}$&2.9814$\times 10^{-9}$&$\sim$4.0882 \\ 
										&&&$u_x$&2.5774$\times 10^{-6}$&1.2644$\times 10^{-7}$&7.3962$\times 10^{-9}$&4.5508$\times 10^{-10}$&$\sim$4.1558\\ 
										&&&$u_y$&1.3175$\times 10^{-6}$&6.9426$\times 10^{-8}$&4.1467$\times 10^{-9}$&2.5654$\times 10^{-10}$&$\sim$4.1088 \\ 
										&&&$u_z$&2.1628$\times 10^{-6}$&9.9075$\times 10^{-8}$&5.5962$\times 10^{-9}$&3.4079$\times 10^{-10}$&$\sim$4.2106 \\ 
										\cline{4-9}
										\multirow{4}{*}{0.10}&\multirow{4}{*}{1/10}&\multirow{4}{*}{1/25}&$\theta$&1.3769$\times 10^{-7}$&7.9947$\times 10^{-9}$	&4.9629$\times 10^{-10}$&3.1115$\times 10^{-11}$&$\sim$4.0372\\ 
										&&&$u_x$&4.5349$\times 10^{-8}$&7.9947$\times 10^{-9}$	&4.9629$\times 10^{-10}$&3.1115$\times 10^{-11}$&$\sim$4.0727\\ 
										&&&$u_y$&2.1041$\times 10^{-8}$&1.1906$\times 10^{-9}$&7.2976$\times 10^{-11}$&4.5464$\times 10^{-12}$&$\sim$4.0587 \\ 
										&&&$u_z$&5.8170$\times 10^{-8}$&3.3611$\times 10^{-9}$&	2.0675$\times 10^{-10}$&
										1.2891$\times 10^{-11}$&$\sim$4.0466 \\ 
										\cline{4-9}
										\multirow{4}{*}{0.15}&\multirow{4}{*}{1/10}&\multirow{4}{*}{1/25}&$\theta$&1.4106$\times 10^{-9}$&
										8.1745$\times 10^{-11}$&5.0850$\times 10^{-12}$&3.3828$\times 10^{-13}$&$\sim$4.0086\\
										&&&$u_x$&6.5890$\times 10^{-10}$&3.7145$\times 10^{-11}$&2.2804$\times 10^{-12}$&1.4216 $\times 10^{-13}$&$\sim$4.0594\\
										&&&$u_y$&2.5730$\times 10^{-10}$&1.4675$\times 10^{-11}$&9.0264$\times 10^{-13}$&5.6299$\times 10^{-14}$&$\sim$4.0527\\ 
										&&&$u_z$&9.7287$\times 10^{-10}$&5.6602$\times 10^{-11}$&2.1824$\times 10^{-13}$&3.4959$\times 10^{-12}$&$\sim$4.0407   \\ 
										\cline{4-9}
										\multirow{4}{*}{0.20}&\multirow{4}{*}{1/10}&\multirow{4}{*}{1/25}&$\theta$&6.4539$\times 10^{-11}$&
										2.9193$\times 10^{-12}$&1.7533$\times 10^{-13}$&9.5754$\times 10^{-15}$&$\sim$4.2395\\
										&&&$u_x$&3.8868$\times 10^{-11}$&1.7540$\times 10^{-12}$&1.0422$\times 10^{-13}$&6.5206$\times 10^{-15}$&$\sim$4.1804\\ 
										&&&$u_y$&1.7096$\times 10^{-11}$&	7.9028$\times 10^{-13}$&4.7037$\times 10^{-14}$&3.0432$\times 10^{-15}$&$\sim$4.1519 \\ 
										&&&$u_z$&7.4620$\times 10^{-11}$&3.3650$\times 10^{-12}$&1.9973$\times 10^{-13}$&1.2351$\times 10^{-14}$&$\sim$4.1869 \\ 
										\hline\hline
									\end{tabular}
								\end{center} 
							\end{table} 
							\noindent\textbf{Example 4.} We further consider the four-dimensional case that can be used to test the generality of the developed fourth-order MRT-LB model for the $d$-dimensional  coupled Burgers' equations, where   the periodic condition is considered and the following initial conditions are given by  
							\begin{align}
								\begin{cases}
									u_{x_1}(x_1,x_2,x_3,x_4,0)=-2\upsilon\frac{   \pi \cos( \pi x_1)\sin(2\pi x_2)\sin(3\pi x_3)\sin(4\pi x_4) }{2+\sin(\pi x_1)\sin(2\pi x_2)\sin(3\pi x_3)\sin(4\pi x_4) },&(x_1,x_2,x_3,x_4)\in\Omega,\\
									u_{x_2}(x_1,x_2,x_3,x_4,0)=-2\upsilon\frac{  2\pi \sin(\pi x_1)\cos(2\pi x_2)\sin(3\pi x_3)\sin(4\pi x_4) }{2+\sin(\pi x_1)\sin(2\pi x_2)\sin(3\pi x_3)\sin(4\pi x_4) },&(x_1,x_2,x_3,x_4)\in\Omega,\\
									u_{x_3}(x_1,x_2,x_3,x_4,0)=-2\upsilon\frac{  3\pi \sin(\pi x_1)\sin(2\pi x_2)\cos(3\pi x_3)\sin(4\pi x_4) }{2+\sin(\pi x_1)\sin(2\pi x_2)\sin(3\pi x_3)\sin(4\pi x_4) },&(x_1,x_2,x_3,x_4)\in\Omega,\\
									u_{x_4}(x_1,x_2,x_3,x_4,0)=-2\upsilon\frac{  4\pi \sin(\pi x_1)\sin(2\pi x_2)\sin(3\pi x_3)\cos(4\pi x_4) }{2+\sin(\pi x_1)\sin(2\pi x_2)\sin(3\pi x_3)\sin(4\pi x_4) },&(x_1,x_2,x_3,x_4)\in\Omega,
								\end{cases}
							\end{align}
							where the computational domain is $\Omega=\{(x_1,x_2,x_3,x_4):-1\leq x_1,x_2,x_3,x_4\leq 1\}$. The analytical solution of $\bm{u}=(u_{x_1},u_{x_2},u_{x_3},u_{x_4})^T$ can be obtained as
							\begin{align}
								\left\{\begin{aligned}
									&u_{x_1}=-2\upsilon\frac{  \pi \exp{\big(-30\upsilon \pi^2 t\big)}\cos( \pi x_1)\sin(2\pi x_2)\sin(3\pi x_3)\sin(4\pi x_4)}{2+\exp{\big(-30\upsilon \pi^2 t\big)}\sin( \pi x_1)\sin(2\pi x_2)\sin(3\pi z)\sin(4\pi x_4)},\\
									&u_{x_2}=-2\upsilon\frac{ 2\pi \exp{\big(-30\upsilon \pi^2 t\big)}\sin( \pi x_1)\cos(2\pi x_2)\sin(3\pi x_3)\sin(4\pi x_4)}{2+\exp{\big(-30\upsilon \pi^2 t\big)}\sin(2\pi x_1)\sin(4\pi x_2)\sin(6\pi z)\sin(8\pi x_4)},\\
									&u_{x_3}=-2\upsilon\frac{ 3\pi \exp{\big(-30\upsilon \pi^2 t\big)}\sin( \pi x_1)\sin(2\pi x_2)\cos(3\pi x_3)\sin(4\pi x_4)}{2+\exp{\big(-30\upsilon \pi^2 t\big)}\sin( \pi x_1)\sin(2\pi x_2)\sin(3\pi z)\sin(4\pi x_4)},\\
									&u_{x_4}=-2\upsilon\frac{ 4\pi \exp{\big(-30\upsilon \pi^2 t\big)}\sin( \pi x_1)\cos(2\pi x_2)\sin(3\pi x_3)\cos(4\pi x_4)}{2+\exp{\big(-30\upsilon \pi^2 t\big)}\sin( \pi x_1)\sin(2\pi x_2)\sin(3\pi z)\sin(4\pi x_4)},
								\end{aligned}\right.
							\end{align}
							from which one can obtain the analytical solution of $\theta(x_1,x_2,x_3,x_4,t)$,
							\begin{align}
								\theta(x_1,x_2,x_3,x_4,t)=2+\exp\big(-30\upsilon\pi^2t\big)\sin(\pi x_1)\sin(2\pi x_2)\sin(3\pi x_3)\sin(4\pi x_4).
							\end{align}
							Here we take different values of parameter $\varepsilon$ ($\varepsilon=$0.08, 0.10)  and time $t$ ($t=$1.0, 2.0)   to calculate the RMSEs and CRs, where the lattice spacing is varied from $\Delta x=1/10$ to $\Delta x=1/25$ with $\Delta x^2/\Delta t=0.25$. As shown in Table \ref{Table-4D-1}, the proposed MRT-LB model has a fourth-order convergence rate for the four-dimensional coupled Burgers' equations.
							\begin{table} 
								\begin{center}
									\caption{The RMSEs and CRs under different values of parameter  $\varepsilon$  and time $t$ ($\Delta x=1/10$, $\Delta t=1/40$).} 
									\label{Table-4D-1}
									\begin{tabular}{cccccccccccccc}\hline\hline
										$\varepsilon$&$t$&variable&RMSE$_{\Delta x,\Delta t}$&RMSE$_{2\Delta x/3,4\Delta t/9}$&RMSE$_{\Delta x/2,\Delta t/4}$&RMSE$_{2\Delta x/5,4\Delta t/25}$&CR\\
										\hline
										\multirow{4}{*}{0.08}&\multirow{4}{*}{1.0}&$\theta$& 3.7475$\times 10^{-7}$&6.0854$\times 10^{-8}$&2.1042$\times10^{-8}$ &8.5486$\times 10^{-9}$&$\sim$ 4.1329 \\ 
										&&$u_{x_1}$&6.7239$\times 10^{-8}$&1.2084$\times 10^{-8}$& 3.7095$\times 10^{-9}$&1.5000$\times 10^{-9}$&      $\sim$ 4.1569   \\ 
										&&$u_{x_2}$&1.3460$\times 10^{-7}$  &2.3975$\times 10^{-8}$&7.3399$\times 10^{-9}$ &2.9646$\times 10^{-9}$&$\sim$ 4.1715 \\ 
										&&$u_{x_3}$&2.0032$\times 10^{-7}$ &3.4553$\times 10^{-8}$&1.0466$\times 10^{-8}$&4.2069$\times 10^{-9}$&$\sim$ 4.2258 \\ 
										&&$u_{x_4}$&2.1478$\times 10^{-7}$ &3.7529$\times 10^{-8}$&1.1409$\times 10^{-8}$&4.5931$\times 10^{-9}$&$\sim$ 4.2049\\ 
										\cline{3-8}
										\multirow{4}{*}{0.10}&\multirow{4}{*}{1.0}&$\theta$&2.7710$\times 10^{-8}$&5.1002$\times 10^{-9}$&1.5861$\times 10^{-9}$&6.4622$\times 10^{-10}$&$\sim$ 4.1083 \\ 
										&&$u_{x_1}$&4.7313$\times 10^{-9}$&8.5449$\times 10^{-9}$&2.6283$\times 10^{-10}$&1.0639$\times 10^{-10}$&$\sim$ 4.1480 \\ 
										&&$u_{x_2}$&1.0345$\times 10^{-8}$&1.8709$\times 10^{-9}$&5.7627$\times 10^{-10}$&2.3345$\times 10^{-10}$&$\sim$ 4.1443 \\ 
										&&$u_{x_3}$&1.6755$\times 10^{-8}$&2.9399$\times 10^{-9}$&8.9723$\times 10^{-10}$&3.6201$\times 10^{-10}$&$\sim$ 4.1942 \\ 
										&&$u_{x_4}$&1.6801$\times 10^{-8}$&3.0219$\times 10^{-9}$&9.2956$\times 10^{-10}$&3.7638$\times 10^{-10}$&$\sim$ 4.1529 \\   
										\cline{3-8}
										\multirow{4}{*}{0.08}&\multirow{4}{*}{2.0}&$\theta$&5.7439$\times 10^{-11}$ &1.0325$\times 10^{-11}$& 3.1931$\times 10^{-12}$ &1.2982$\times 10^{-12}$&$\sim$4.1450      \\ 
										&&$u_{x_1}$&1.1171$\times 10^{-11}$&1.9856$\times 10^{-12}$&6.0973$\times 10^{-13}$&2.4676$\times 10^{-13}$&$\sim$4.1698   \\ 
										&&$u_{x_2}$&2.2350$\times 10^{-11}$  &3.9563$\times 10^{-12}$&1.2134$\times 10^{-12}$&4.9080$\times 10^{-13}$&$\sim$4.1767 \\ 
										&&$u_{x_3}$&3.3400$\times 10^{-11}$&5.8255$\times 10^{-12}$& 1.7783$\times 10^{-12}$&7.1776$\times 10^{-13}$&$\sim$4.2013 \\ 
										&&$u_{x_4}$&4.0353$\times 10^{-11}$&7.1068$\times 10^{-12}$&2.1752$\times 10^{-12}$&8.7895$\times 10^{-13}$&$\sim$4.1857  \\ 
										\cline{3-8}
										\multirow{4}{*}{0.10}&\multirow{4}{*}{2.0}&$\theta$& 3.7032$\times 10^{-13}$ &6.7464$\times 10^{-14}$&		 
										2.1124$\times 10^{-14}$  &8.7585$\times 10^{-15}$&$\sim$4.0926  \\ 
										&&$u_{x_1}$&7.8597$\times 10^{-14}$ &1.4108$\times 10^{-14}$&4.3447$\times 10^{-15}$&1.8801$\times 10^{-15}$&$\sim$4.0556  \\ 
										&&$u_{x_2}$&1.6371$\times 10^{-13}$&2.9383$\times 10^{-14}$& 9.0502$\times 10^{-15}$&3.7256$\times 10^{-15}$&$\sim$4.1303  \\ 
										&&$u_{x_3}$&2.5471$\times 10^{-13}$ &4.5044$\times 10^{-14}$& 1.3810$\times 10^{-14}$ &5.6307$\times 10^{-15}$&$\sim$4.1662  \\ 
										&&$u_{x_4}$&2.9869$\times 10^{-13}$ &5.3572$\times 10^{-14}$&1.6493$\times 10^{-14}$&6.7207$\times 10^{-15}$&$\sim$ 4.1466\\ 
										\hline\hline
									\end{tabular}
								\end{center} 
							\end{table}

							\section{Conclusions}\label{Conclusion}	
							In this paper, we first adopted the Cole-Hopf transformation to eliminate the nonlinear convection terms in  the $d$-dimensional  ($d\geq 1$) coupled Burgers' equations, and in this case, a simple diffusion equation is obtained. In particular, the velocity $\bm{u}$ in the $d$-dimensional coupled Burgers' equations can be determined by the variable $\theta$ and  gradient term $\nabla\theta$ in the transformed diffusion equation. Then at the diffusive scaling,  we obtained the macroscopic finite-difference scheme   of the MRT-LB model for the $d$-dimensional transformed diffusion equation, and also the consistent fourth-order modified equation  of the MRT-LB model  through the Maxwell iteration method. Furthermore, under the fourth-order approximation of the distribution function $\bm{f}$, we obtained the fourth-order MRT-LB model for the $d$-dimensional coupled Burgers' equations. It is worth mentioning that under the condition (\ref{Bur-Sol}), the gradient term  $\nabla\theta$ in the $d$-dimensional  diffusion equation  can also be calculated by the non-equilibrium distribution function  with a fourth-order accuracy. What is more, it should be noted that the developed fourth-order MRT-LB model for the coupled Burgers' equations is not limited to the spatial dimension, which indicates the universality of the LB method. Finally, we carried out some numerical experiments  of the coupled Burgers' equations  to test the developed MRT-LB model, and the numerical results show  that the proposed MRT-LB model has a fourth-order convergence rate in space, which is also consistent with the theoretical analysis.
							\section*{Acknowledgements}
							The computation is completed in the HPC Platform of Huazhong University of Science and Technology. This work was financially supported by the National Natural Science Foundation of China (Grants No. 12072127 and No. 51836003), Interdiciplinary
							Research Program of Hust (2023JCJY002) and the Fundamental Research Funds for the Central Universities, Hust (No. 2023JYCXJJ046).

						\end{document}